\documentclass[a4paper,12pt]{article}
\usepackage[dvips]{epsfig}
\usepackage{amsmath,amssymb,euscript,graphicx,amsfonts,pifont}
\usepackage{enumerate,color}
\usepackage{graphicx}
\usepackage{hyperref}
\usepackage{tikz}
\usetikzlibrary{calc, decorations.markings, arrows}
\newtheorem{theorem}{Theorem}[section]
\newtheorem{proposition}[theorem]{Proposition}
\newtheorem{lem}[theorem]{Lemma}

\newcounter{case}

\renewcommand{\thecase}{\arabic{case}}

\newcounter{subcase}

\numberwithin{subcase}{case}

\def\a{\alpha} \def\b{\beta} \def\g{\gamma} \def\d{\delta} 
  \def\t{\tau}  
\def\th{\theta} \def\ld{\lambda} \def\ph{\varphi}
\def\Del{\Delta}  
 \def\O{\Omega} \def\G{\Gamma}
\def\D{{\rm \Delta }}

 \def\ou{\overline U}

\def\o{\overline}  
\def\di{\bigm|} \def\lg{\langle} \def\rg{\rangle}

\def\PSL{\hbox{\rm PSL}}  
\def\Aut{\hbox{\rm Aut\,}}  
 \def\soc{\hbox{\rm soc}}

 \def\PG{\hbox{\rm PG}} \def\PGL{\hbox{\rm PGL}}
  \def\GL{\hbox{\rm GL}}  \def\P\GL{\hbox{\rm P\GL}}
  \def\FF{{\hbox{\sf F\kern-.43emF}}}
\def\cal{\mathcal}  
\def\PGammaL{{\rm P\Gamma L}}  \def\PSp{{\rm PSp}}
\def\Sp{{\rm Sp}} \def\PSU{{\rm PSU}} \def\SU{{\rm SU}}  \def\GU{{\rm GU}} 
 \def\diag{{\rm diag}}

\def\orho{\overline{\rho}}

\def\P{\hbox{\rm P}}

\def\nd{\mathrel{\bigm|\kern-.7em/}} 
 \def\f{\noindent}
\def\qed{\hfill $\Box$} \def\demo{\f {\bf Proof}\hskip10pt}

\begin{document}

	\begin{center}
		{\bf\large  Hamilton Cycles In
			Graphs  of Rank 3

 (classical group cases (I))}
		\footnote{This work was supported in part  by the National Natural Science Foundation of China
			(12471332).}
	\end{center}
	
	%\bigskip
	
	\begin{center}
		Shaofei Du and Chen Yao\footnote{Corresponding author: 2240501007@cnu.edu.cn}\\
		\medskip
		{\it {\small
				School of Mathematical Sciences, Capital Normal University,\\
				Bejing 100048, People's Republic of China\\
		}}
	\end{center}

	\renewcommand{\thefootnote}{\empty}%{footnote}}

\begin{abstract}   This paper contributes to a long standing Lov\'asz  problem regarding existence of Hamilton paths (Hamilton cycles, asked by Alspach) in vertex-transitive graphs.
 A connected simple graph is said to be of rank $r$ if  its automorphism group is a transitive group of   rank $r$ on vertices.
 Since  every  graph of rank 2  is complete,  a basic step is to  determine  the Hamiltonicity of graphs  of  rank 3.
  It is easy to see that every   graph  of rank 3  whose  automorphism group  is imprimitive must contain a  Hamilton cycle.
 In this paper,   it will be shown that with the exception  of the Petersen graph,  every connected graph of rank 3 contains a  Hamilton cycle,
  provided its automorphism group $G$  is  a  primitive group of rank 3,  whose socle $\soc(G)$ is a linear group,    symplectic group or unitary   group.
 The remaining  primitive  groups of rank 3  will be discussed   in our other papers.
\end{abstract}

\section{Introduction}
\label{sec:intro}
\indent
A simple path (resp. cycle) going through all vertices of a graph is called a Hamilton path (resp. cycle).
In 1969, Lov\'asz \cite{L70} asked that if there exists  a
finite  connected vertex-transitive graph of order more than 3 without a Hamilton path; and   in 1981, Alspach  \cite{A81} asked if there exists an infinite number of connected vertex-transitive graphs that do not have a Hamiltonian cycle.
Till now,  with the exception of $K_2$, only four connected
vertex-transitive graphs that do not have a Hamilton cycle are known to exist.
These four graphs are the Petersen graph, the Coxeter graph and the two graphs obtained
from them by replacing each vertex by a triangle.
The fact that none of these four graphs
is a Cayley graph  has led to a folklore conjecture that except for $K_2$,
every Cayley graph has a Hamilton cycle.

The existence of Hamilton paths, and in some cases
also Hamilton cycles, in  connected vertex-transitive graphs
has been shown for graphs of particular orders, such as, $kp$, where $k\le 6$,
$p^j$, where $j\le 5$, and $2p^2$, where throughout this paper $p$  denotes a prime, see \cite{C98,DKM21,DTY23,KM08,KMZ12,KS09,DM87,DM92,DM85,MP82,MP83,MS4,MS2,MS5,Z15}
and a survey paper \cite{KM09}.  For known results of  Hamiltonicity  of  Cayley graphs, see    \cite{A82,CG96,D83,
	GWM14,GM07,GKMM12,KM09,MR234, DM83,WM18,WM15}. Additionally, the existence of Hamilton cycle of semisymmetric graphs of order $2pq$ and $2p^3$ were recently confirmed in \cite{Chen26,DY26}, where $p$ and $q$ are primes.

 A connected simple graph is said to be of rank $r$ if  its automorphism group is a transitive group of   rank $r$ on vertices.
     Every   graph of rank 2  is a complete graph,   which contains a Hamilton cycle, and  so graphs of rank 3 are the smallest nontrivial case.
  One of our long-term goals is to show that every  such graph   contains a Hamilton cycle, except for the Petersen graph.
  It is straight to see  the the existence of   Hamilton cycles when
 the group is imprimitive. So we only need to deal with primitive groups of rank 3.
  In this paper,  we shall deal with primitive groups of rank 3, whose socle is a  linear group,    symplectic group or unitary   group.
   The other remaining cases will be done in future.
The following  theorem is the main result of  this paper.

\vskip -2mm
\begin{theorem} \label{the:main} Let $\G$ be a graph of rank 3, whose automorphism group   is a  primitive group of rank 3, whose socle is a linear group,    symplectic group or unitary   group.
 Then  $\G $ contains a Hamilton cycle, except for  the Petersen graph.
\end{theorem}

%\begin{theorem} \label{the:main2} Let $\G$ be a graph of rank 3, whose automorphism group  $G$  is one of the following classes: (1) imprimitive groups;  (2) affine primitive groups; (3) $\soc(G)=T\times T$, where $T$ is simple; (4) $\soc(G)=A_n$. Then  $\G $ contains a %Hamilton cycle, except for  the Petersen graph.
%\end{theorem}
\vskip -2mm
Remak that the existence of Hamilton cycles  has been confirmed   for  several graphs of rank 3, for instance,
  graphs of order $pq$ and rank 3 where $p$ and $q$ are primes, see \cite{DKM21}.

The notation in this paper are standard, see \cite{BHD67,Hup67}.
This paper is organized as follows: after this introductory section, some preliminaries     will be given in Section 2; the strategy of the proof of Theorem~\ref{the:main}  will be given in Section 3;  and
the   primitive  groups of rank 3, whose socle is a  linear group, symplectic group or  unitary
 will be discussed  in Sections 4-6, separately.
\section{Preliminaries}
\vskip -2mm
 \label{sec:pre}
Throughout this paper, graphs are finite, simple, connected and undirected,
and groups are finite.
Furthermore, a {\em multigraph} is a generalization
of a graph in which multiple edges and loops are allowed.
Given a graph $X$, we denote the vertex set by $V(X)$ and the edge set by $E(X)$, respectively.
For two adjacent vertices $u,v \in V(X)$, we write
$u \sim v$ and denote the corresponding edge by $\{u,v\}$.

 A permutation group $G$ on a set $V$ induces
the action of $G$ on $V\times V$, whose
orbits are called {\em orbitals} of $G$, while  $O_0=\{(x,x) \colon x\in V\}$ is said to be {\it trivial}.
Every orbital $O_i$ corresponds to a {\it paired} orbital $O_i^*=\{ (x,y): (y,x)\in O_i\}$, and
$O_i$ is said to be {\it self-paired} if $O_i=O_i^*$.

For an arbitrary union $\cal{O}$ of orbitals (having empty intersection
with $O_0$), the {\em generalized orbital (di)graph} $X(V,\cal{O})$
of the action of $G$ on $V$
with respect to $\cal{O}$ is a simple (di)graph with the vertex set $V$ and
edge set $\cal{O}$.
It can be viewed as an undirected graph if and only if  $\cal{O}$ coincides with its symmetric closure,
that is, $\cal{O}$ has the property that $(x,y)\in \cal{O}$ implies $(y,x)\in \cal{O}$.
Further,  $X(V,\cal{O})$ is
said to be {\em basic orbital graph } if ${\cal O} =O_i\cup O_i^*$ for some $O_i$.

 Let $X$ be a graph with a transitive group $G$ of $\Aut(X)$. By $G_v$ and $[G:H]$ we denote the stabilizer of $v\in V(X)$ and the set of cosets of $G$ relative to $H$, respectively, where $H=G_v$.
Then there exists a 1-1 correspondence between orbitals and suborbits of $G$, that is orbits of $H$ on $V(X)$. A suborbit corresponding to a {\em self-paired}
orbital is said to be self-paired. Therefore, the (generalized) orbital (di)graph $X(V, \cal{O})$ is denoted by  $X(G,H,D)$, provided $\cal{O}$ corresponds to a union   of suborbits  (some double cosets of $H$).
 \vskip 1.5mm
Finally, we give some  known results on graphs and groups.

\begin{proposition}
	\label{pro:1}
	{\rm \cite[Theorem~1]{BJ78}}
	Every $2$-connected regular graph of order $n$ and valency at least $n/3$ contains a Hamilton cycle (every vertex-transitive graph is 2-connected, see \cite[Theorem~3]{M70}).
\end{proposition}

\begin{proposition}
    \label{pro:2}
    {\rm \cite[Corollary~3]{A82}}
    If $X$ is a connected Cayley graph of an abelian group of order at least 3, then every edge of $X$ lies in a Hamiltonian cycle.
\end{proposition}

\begin{proposition} \label{pro:3}  {\rm \cite{KL82}} Let $\soc(F/Z(M))=M/Z(M)$ be a linear group,    symplectic group or unitary   group,  acting as a primitive
rank $3$ permutation group of degree $n$, on the set $X$ of cosets of a subgroup $K/Z(M)$ of $F/Z(M)$. Then
one of the following occurs up to conjugacy under $\Aut (M/Z(M))$:
\begin{enumerate}
\item[\rm(1)] $M/Z(M)=\PSL(n,q)\le F/Z(M)\le \Aut(M/Z(M))$ and $X$ is the set of lines, $n\ge 4;$
\item[\rm(2)]  $M$ is one of the groups
$\Sp(2m-2,q)$, $\SU(m,q)$,
for $m\ge 3$ and q a prime power,  $M\trianglelefteq F$ with $F/Z(M)\le\Aut(M/Z(M));$
\begin{enumerate}
\item[\rm(2.1)]  $X$ is an $M$-orbit of singular (or  isotropic) points;
\item[\rm(2.2)] $X$ is an $M$-orbit of maximal totally singular (or  isotropic) subspaces and
$M=\Sp(4, q), $ $\SU(4, q), \SU(5, q).$
\item[\rm(2.3)] $X$ is any $M$-orbit of nonsingular points and $M=\SU(m,2).$
\end{enumerate}
\item[\rm(3)] groups of small order

$M/Z(M)=\PSL(2,4)$, $|X|=10$; $M/Z(M)=\PSL(2,9)\cong A_6,$ $|X|=15$;

$M/Z(M)=\PSL(4,2)\cong A_8,$ $|X|=28$; $F/Z(M)=\PGammaL(2,8)$ $|X|=36$;

$M/Z(M)=\PSL(3,4)$, $M/Z(M)\cap (K/Z(M))=A_6$;

$M/Z(M)=\PSL(4,3)$, $M/Z(M)\cap (K/Z(M))=\PSp(4,3)$;

$M=\Sp(6,2), K=G_2(2);$  $M=\SU(3,3),$ $K\cap  M=\PSL(3,2);$

$M =\SU(3,5), K\cap M = 3\cdot A_7$; $M =\SU(4,3)$, $K\cap M=4\cdot \PSL(3,4)$;

$M=\SU(6,2), K\cap M=3\cdot \PSU(4,3)\cdot 2.$
\end{enumerate}
\end{proposition}

\section{Proof of   Theorem~\ref{the:main}}
Before proving Theorem\ref{the:main}, we first deal with graphs of rank 3, with  the imprimitive  automorphism group.
\begin{lem} Suppose that $\G$ is a connected graph of rank 3,   with the imprimitive  automorphism group $G=\Aut(\G)$. Then $\G$ contains a Hamilton cycle.
 \end{lem}
\demo  Suppose that  $G$  is imprimitive. Let $\bf B$ be a complete block system. Take a point $\a$ in a block $B$. Since $G$ is of rank 3,    two nontrivial   orbits of $G_\a$ are
  $\Del_1=\b^{G_\a}$ for any $\b\in B\setminus \{\a\}$ and  $\Del_2=V\setminus B$, with respective length  $|B|-1$ and $|V|-|B|$. Clearly, $|\Del_1|\lvertneqq \frac{|V|}2$ and $|\Del_2|\ge   \frac{|V|}2$.
Then the corresponding basic orbital graph $\G(V, \Del_1)$ is disconnected.
 By Proposition~\ref{pro:1}, $\G(V, \Del_2)$ has a Hamilton cycle.     \qed

\vskip 3mm
\f {\bf Proof of  Theorem~\ref{the:main}:} To prove Theorem~\ref{the:main}, let $\G$ be a  graph of rank 3, whose automorphism group   $G=\Aut(\G)$ is  a primitive  of rank 3, such that
  $\soc(G)=M/Z(M)$ is  a  linear group,    symplectic group or unitary   group.   Then $M$ is one of groups listed in Proposition~\ref{pro:3}.
All groups  in  (3)   with small order  may be managed by Magma. More efficiently,   first pick up a semiregular automorphism $\lg \rho \rg $, then study the $\lg \rho \rg $-blocks graph by Magma to find a Hamilton cycle for that,  which may be lifted to a  Hamilton cycle of $\G$.
In the following three sections, we shall the theorem  for the groups in  the  (1) and  (2) of  Proposition~\ref{pro:3},  according to  $M$ is a linear group, symplectic group and unitary group, separately.
It will be shown that in all cases, $\G$ contains a  Hamilton cycle. Thus Theorem~\ref{the:main} is proved.
\qed

\vskip 3mm
Throughout the paper,  let $\FF_q^*=\lg \th  \rg$ and  $\FF_{q^n}^*=\lg \ld \rg$, where $n\ge 2$,  $\ld ^n=\sum_{i=0}^{n-1}a_i\ld ^i$. Then $\overline {\ld }(x)=\ld x$,  $\forall x\in \FF_{q^n}$  gives a linear transformation of $\FF_{q^n}$ over $\FF_q$.
Moreover, the matrix $A$ of $\overline{\ld }$  under the base: $1, \ld, \cdots, \ld^{n-1}$ is
\begin{equation}\label{A1} A={\left(\begin{array}{cccccc}0 &1 &0 &\cdots &0 &0\\ 0 &0 &1 &\cdots &0 & 0 \\  &&&\cdots&&\\0 &0 &0  &\cdots &0&1 \\   a_0 &a_1 &a_2 &\cdots
&a_{n-2} &a_{n-1} \end{array}\right)}\end{equation}
and then
\begin{equation}\label{A2} A^{-1}={\left(\begin{array}{ccccc} -a_1a_0^{-1} &-a_2a_0^{-1} &\cdots &-a_{n-1}a_0^{-1} & a_0^{-1}\\ 1 &0 &\cdots &0 &0\\ 0 &1 &\cdots &0 &0\\ &&\cdots &&\\0 &0 &\cdots &1 &0\end{array}\right)}.\end{equation}

 By $\Del_s$ and $\Del_l$ we denote the shorter and longer one, respectively, where $|\Del_s|$  may be equal to  $|\Del_l|.$ Then  $|\Del_l|\ge \frac{|V|}3.$
  By Proposition~\ref{pro:1}  the corresponding basic orbital graph $\G(V, \Del_l)$ has a Hamilton cycle.  So we one need to prove the theorem for  $\G(V, \Del_s)$. Later we always denote $V$ by $\O$. So
 set
\begin{equation} \label{orbital} \G=\G(\O, \Del_s) .  \end{equation}

 For any H-path leading $\a$ to $\b$  of a certain graph,  we can simply denote that by $\a \to \b$. From now on,   use the abbreviation   H-cycle (resp. H-path) for  Hamilton cycle (resp. Hamilton path).
Finally  by $A^\t$  we denote   the inverse transpose of $A$. Address that we always identify  a subspace $W$ in $V=V(n,q)$ with  the corresponding subspace  $\PG(W)$ of the projective space $\PG(V)$ of $V$.

\section{$M$ is   a linear group}
\vskip -2mm
By Proposition~\ref{pro:3}, we shall deal with a action for the linear groups in the following subsection.

 \begin{lem}\label{L1} Let $\soc(G)=\PSL(n, q)$, where $n\ge 4$, acting on the set $\O$  of lines.
Then  $\G(\O, \Del_s)$  has a H-cycle.
\end{lem}
\demo  Take a line $\ell $ in $\O$, where   $|\O|=\frac{(q^n-1)(q^{n-1}-1)}{(q^2-1)(q-1)}$.
Then  for any group $G$ such that $M/Z(M)=\PSL(n, q)\le G\le \Aut(M/Z(M))$,  it  is of  rank 3.
Pick up a line $\ell$ and let $\D_l$ and $\D_s$ be the related longer and shorter nontrivial  suborbits, respectively.  Let $\G=\G(\O, \Del_s)$, see \ref{orbital}.

If $(n,q)=(4.2)$, then  $|\O|=35$ and $\D_s=\{ \ell'\in \O\di d(\ell\cap \ell')=0\}$ with size 16.
 By Proposition~\ref{pro:1},  $\G$  has a H-cycle.
So suppose  $(n, q)\ne (4,2)$ so that  $\Del_s=\{ \ell'\in \O\di d(\ell\cap \ell')=1\}.$  Then
 $\{\ell_1, \ell_2\}\in E(\G)$ if and only if  $d(\ell_1\cap \ell_2)=1$.

Let $\FF_{q^n}^*=\lg \ld \rg $,  $\rho=A$, defined in Eq(\ref{A1}) and $\orho =\rho Z$, where $Z=Z(\GL(n,q))$.  Then $\lg \orho \rg $ acts regularly on all points. Set $r=\frac{q^n-1}{q-1}$.
 We shall deal with two cases when $n$ is even and odd, separately.

\vskip 1.5mm
(1) Suppose $n$ is even.   Then  $q+1\di r$. If a non-trivial subgroup $N$ of $\lg \orho \rg $ fixes a line, then it is isomorphic to a subgroup  of $\PGL(2,q)$ and so $|N|\di (q+1)$.  Therefore  $N\le \lg \a\rg $ where  $\a=\ld ^{\frac r{q+1}}.$ Note $\FF_{q^2}=\lg \a\rg  \cup \{0\}$ is just a subspace of dimension 2 over $\FF_q$.
Suppose that $N$ fixes a line $\ell_0$ containing $\ld ^i$ for some $0 \le i \le q^n-1$. Then $\ell_0=\FF_q\lg \ld^i, \a^j \ld^i \rg= \ld^i \FF_q\lg 1, \a^j \rg = \ld^i \FF_q\lg 1, \a \rg$ for some $ 1 \le j < q+1$.
 Since $\FF_q\lg 1, \a \rg=\FF_{q^2}$,   the number of fixed lines of $N$ is equal to
 the number of cosets of $\lg \a \rg$ in $\FF_{q^n}^*$,  that is   $\frac r{q+1}$.

Let  $$\ell_0=\PG(\FF_{q^2})=\{\lg \g\rg \di \g\in \FF_{q^2}\}\,\, {\rm and}\,\,  B_0=\{ \ell_0^{\orho^i}\di 0\le i \le \frac r{q+1}-1\}.$$ Then  $\lg \orho\rg $ acts transitively on the points  in $B_0$, with the stabilizer $\lg \orho^{\frac{r}{q+1}}\rg$, and
  acts  semiregularly on other lines, while  let ${B}_1, {B}_2, \cdots , {B}_k$  be  all  orbits of $\lg \rho\rg $ on $\O\setminus B_0$, where
 $\frac r{q+1}+kr=|\O|$, that is $k=q\cdot \frac {q^{n-2}-1}{q^2-1}$. Therefore,
 $$\O =B_0\cup B_1\cup B_2\cdots \cup B_k.$$

  Take a point $\a$ in $\ell_0=\PG(\FF_{q^2})\in B_0$.   Set $\ell_1=\lg \a, \a^{\orho^{-1}}\rg $.  Since $\ld \not\in \FF_{q^2}$, we have $\a^{\orho^{-1}}\not\in  \FF_{q^2}$ so that  $\ell _1\not\in B_0$, where  no loss, set $\ell_1\in  B_1$.
   In fact, since $\lg \orho \rg$ is transitive on all points of $\PG(n-1, q)$,  every orbit  $B_i$    contains  some lines passing  $\a$, say $\ell_i$, while $\ell_0$ and $\ell_1$ are given as above.
  Then  $B_i=\{ \ell_i^{\orho ^j}\di 0\le j\le r-1\}$.

   Let $\G'$   be the induced subgraph $\G(\cup_{1\le i\le k} B_i)$.  We shall first find a H-cycle $C'$ for  $\G'$ and then  construct a H-cycle $C$  for $\G$ by inserting vertices of $B_0$ to the cycle $C'$.
  Since $\a$ is contained in every $\ell_i$  and  the line $\ell_k$ is adjacent to $\ell_1^{\orho}$, we have  the following H-cycle $C'$ for $\G'$, see Figure 1:

   $$\ell_1, \ell_2, \cdots , \ell_k; \ell_1^{\orho}, \ell_2^{\orho}, \cdots , \ell_k^{\orho};  \ell_1^{\orho^2}, \ell_2^{\orho^2}, \cdots , \ell_k^{\orho^2}; \cdots, \ell_k^{\orho^{r-2}}, \ell_1^{\orho^{r-1}}, \ell_2^{\orho^{r-1}}, \cdots , \ell_k^{\orho^{r-1}}; \ell_1.$$
Since  $\a\in \ell_k\cap \ell_0 \cap \ell_1^{\orho}$, we get    $\ell_k\sim \ell_0\sim \ell_1^{\orho}$ so that  $\ell_k^{\orho^j}\sim \ell_0^{\orho^j}\sim \ell_1^{\orho^{j+1}},$
 where $0\le j \le \frac r{q+1}-1$.  Replacing every edge $\ell_k^{\orho^j},\ell_1^{\orho^{j+1}}$ in the cycle $C'$ by the path $\ell_k^{\orho^j}, \ell_0^{\orho^j}, \ell_1^{\orho^{j+1}},$  we get a H-cycle $C$ for $\G$, see Figure 1.

\vskip 2mm
(2)   Suppose $n$ is odd. Then $(q+1, r)=1$  and so $\lg \orho \rg $ is semiregular on all lines. Therefore,  $B_0$ in (1) is disappeared  and
  $\G$  just coincides with  $\G'$ in (1),  with  the H-cycle $C'$ given by  (1).
\qed

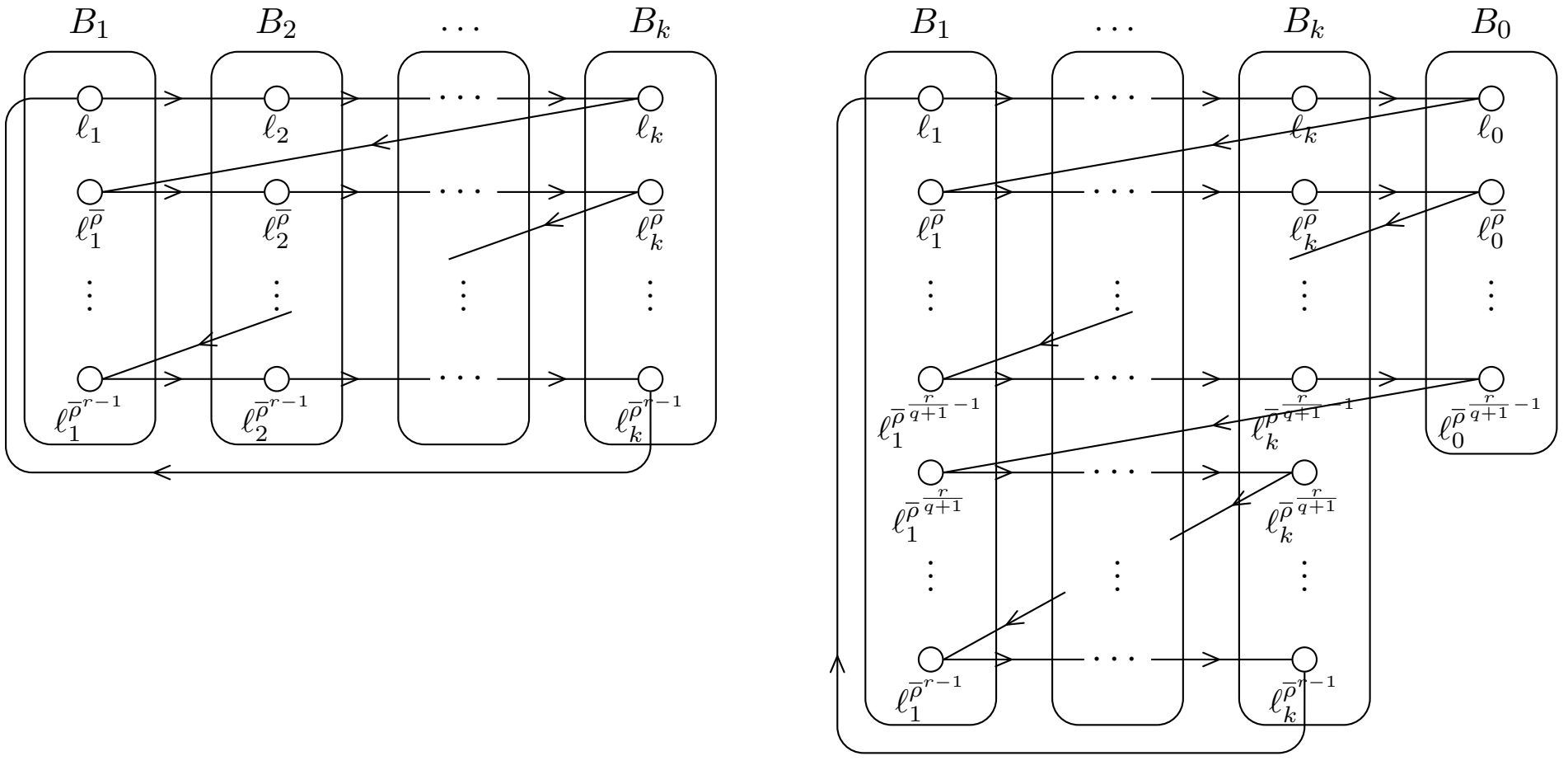
\begin{figure}[htbp]
  \centering
  \begin{tikzpicture}[
    scale=0.75, transform shape,
    p/.style={draw, circle, fill=white, inner sep=0pt, minimum size=2.6mm},
    d/.style={fill=white, inner sep=2pt},
    box/.style={draw, rounded corners=6pt},
    lab/.style={anchor=north, yshift=3pt, font=\footnotesize},
    tlab/.style={anchor=south, font=\small},
    midarrow/.style={
      postaction=decorate,
      decoration={markings, mark=at position 0.5 with {\arrow{angle 45}}}
    },
    remidarrow/.style={
      postaction=decorate,
      decoration={markings, mark=at position 0.5 with {\arrowreversed{angle 45}}}
    }
  ]
    % 111
    \node[p] (a10)  at (0, 0)  {};
    \node[p] (a20)  at (2, 0)  {};
    \node[d] (d10)  at (4, 0)  {$\cdots$};
    \node[p] (a30)  at (6, 0)  {};
    \node[p] (b10)  at (9, 0)  {};
    \node[d] (d20)  at (11, 0)  {$\cdots$};
    \node[p] (b20)  at (13, 0)  {};
    \node[p] (b30)  at (15, 0)  {};
    \draw[midarrow] (a10.east) -- (a20.west);
    \draw[midarrow] (a20.east) -- (d10.west);
    \draw[midarrow] (d10.east) -- (a30.west);
    \draw[midarrow] (b10.east) -- (d20.west);
    \draw[midarrow] (d20.east) -- (b20.west);
    \draw[midarrow] (b20.east) -- (b30.west);
    % 222
    \node[p] (a11)  at (0, -1)  {};
    \node[p] (a21)  at (2, -1)  {};
    \node[d] (d11)  at (4, -1)  {$\cdots$};
    \node[p] (a31)  at (6, -1)  {};
    \node[p] (b11)  at (9, -1)  {};
    \node[d] (d21)  at (11, -1)  {$\cdots$};
    \node[p] (b21)  at (13, -1)  {};
    \node[p] (b31)  at (15, -1)  {};
    \draw[midarrow] (a11.east) -- (a21.west);
    \draw[midarrow] (a21.east) -- (d11.west);
    \draw[midarrow] (d11.east) -- (a31.west);
    \draw[midarrow] (b11.east) -- (d21.west);
    \draw[midarrow] (d21.east) -- (b21.west);
    \draw[midarrow] (b21.east) -- (b31.west);
    % ...
    \node[d] (v1)  at (0, -2)  {$\vdots$};
    \node[d] (v2)  at (2, -2)  {$\vdots$};
    \node[d] (v3)  at (4, -2)  {$\vdots$};
    \node[d] (v4)  at (6, -2)  {$\vdots$};
    \node[d] (v5)  at (9, -2)  {$\vdots$};
    \node[d] (v6)  at (11, -2)  {$\vdots$};
    \node[d] (v7)  at (13, -2)  {$\vdots$};
    \node[d] (v8)  at (15, -2)  {$\vdots$};
    % 333
    \node[p] (a12)  at (0, -3)  {};
    \node[p] (a22)  at (2, -3)  {};
    \node[d] (d12)  at (4, -3)  {$\cdots$};
    \node[p] (a32)  at (6, -3)  {};
    \node[p] (b12)  at (9, -3)  {};
    \node[d] (d22)  at (11, -3)  {$\cdots$};
    \node[p] (b22)  at (13, -3)  {};
    \node[p] (b32)  at (15, -3)  {};
    \draw[midarrow] (a12.east) -- (a22.west);
    \draw[midarrow] (a22.east) -- (d12.west);
    \draw[midarrow] (d12.east) -- (a32.west);
    \draw[midarrow] (b12.east) -- (d22.west);
    \draw[midarrow] (d22.east) -- (b22.west);
    \draw[midarrow] (b22.east) -- (b32.west);

    % 444
    \node[p] (b13)  at (9, -4)  {};
    \node[d] (d23)  at (11, -4)  {$\cdots$};
    \node[p] (b23)  at (13, -4)  {};
    \draw[midarrow] (b13.east) -- (d23.west);
    \draw[midarrow] (d23.east) -- (b23.west);

    % ...
    \node[d] (v11)  at (9, -5)  {$\vdots$};
    \node[d] (v21)  at (11, -5)  {$\vdots$};
    \node[d] (v31)  at (13, -5)  {$\vdots$};
    % 555
    \node[p] (b14)  at (9, -6)  {};
    \node[d] (d24)  at (11, -6)  {$\cdots$};
    \node[p] (b24)  at (13, -6)  {};
    \draw[midarrow] (b14.east) -- (d24.west);
    \draw[midarrow] (d24.east) -- (b24.west);

    % lianjie
    \draw[midarrow] (a30.west) -- (a11.east);
    \draw[midarrow] (a31.west) -- ($(a31)!0.36!(a12)$);
    \draw[midarrow] ($(a31)!0.64!(a12)$) -- (a12.east);
    \draw[midarrow] (b30.west) -- (b11.east);
    \draw[midarrow] (b31.west) -- ($(b31)!0.36!(b12)$);
    \draw[midarrow] ($(b31)!0.64!(b12)$) -- (b12.east);
    \draw[midarrow] (b32.west) -- (b13.east);
    \draw[midarrow] (b23.west) -- ($(b23)!0.36!(b14)$);
    \draw[midarrow] ($(b23)!0.64!(b14)$) -- (b14.east);
    \draw[midarrow, rounded corners=6pt]
    (a32.south) -- (6,-4) -- (-0.9,-4) -- (-0.9,0) -- (a10.west);
    \draw[midarrow, rounded corners=6pt]
    (b24.south) -- (13,-7) -- (8,-7) -- (8,0) -- (b10.west);

    % block
    \draw[box] (-0.7,0.5) rectangle (0.7,-3.7);
    \draw[box] (1.3,0.5) rectangle (2.7,-3.7);
    \draw[box] (3.3,0.5) rectangle (4.7,-3.7);
    \draw[box] (5.3,0.5) rectangle (6.7,-3.7);
    \draw[box] (8.3,0.5) rectangle (9.7,-6.7);
    \draw[box] (10.3,0.5) rectangle (11.7,-6.7);
    \draw[box] (12.3,0.5) rectangle (13.7,-6.7);
    \draw[box] (14.3,0.5) rectangle (15.7,-3.8);

    % sign
    \node[lab] at (a10.south) {$\ell_1$};
    \node[lab] at (a20.south) {$\ell_2$};
    \node[lab] at (a30.south) {$\ell_k$};
    \node[lab] at (a11.south) {$\ell_1^{\orho}$};
    \node[lab] at (a21.south) {$\ell_2^{\orho}$};
    \node[lab] at (a31.south) {$\ell_k^{\orho}$};
    \node[lab] at (a12.south) {$\ell_1^{\orho^{r-1}}$};
    \node[lab] at (a22.south) {$\ell_2^{\orho^{r-1}}$};
    \node[lab] at (a32.south) {$\ell_k^{\orho^{r-1}}$};
    \node[lab] at (b10.south) {$\ell_1$};
    \node[lab] at (b20.south) {$\ell_k$};
    \node[lab] at (b30.south) {$\ell_0$};
    \node[lab] at (b11.south) {$\ell_1^{\orho}$};
    \node[lab] at (b21.south) {$\ell_k^{\orho}$};
    \node[lab] at (b31.south) {$\ell_0^{\orho}$};
    \node[lab] at (b12.south) {$\ell_1^{\orho^{\frac{r}{q+1}-1}}$};
    \node[lab] at (b22.south) {$\ell_k^{\orho^{\frac{r}{q+1}-1}}$};
    \node[lab] at (b32.south) {$\ell_0^{\orho^{\frac{r}{q+1}-1}}$};
    \node[lab] at (b13.south) {$\ell_1^{\orho^{\frac{r}{q+1}}}$};
    \node[lab] at (b23.south) {$\ell_k^{\orho^{\frac{r}{q+1}}}$};
    \node[lab] at (b14.south) {$\ell_1^{\orho^{r-1}}$};
    \node[lab] at (b24.south) {$\ell_k^{\orho^{r-1}}$};

    \node[tlab] at (0,0.5) {$B_1$};
    \node[tlab] at (2,0.5) {$B_2$};
    \node[tlab] at (4,0.5) {$\cdots$};
    \node[tlab] at (6,0.5) {$B_k$};
    \node[tlab] at (9,0.5) {$B_1$};
    \node[tlab] at (11,0.5) {$\cdots$};
    \node[tlab] at (13,0.5) {$B_k$};
    \node[tlab] at (15,0.5) {$B_0$};
  \end{tikzpicture}
  \caption{{\footnotesize  H-cylces  $C'$ of $\G'$ and $C$ of  $\G$ in Lemma 4.1.(1).}}
\end{figure}

\section{$M$ is   a symplectic group}
\vskip -2mm
By Proposition~\ref{pro:3}, we  need to deal with two different  actions  for  the symplectic groups in  the following two  subsections.

\subsection{Action  on isotropic points}
   Let $V(2n,q)$ be the $2n$-dimensional symplectic  space over $\FF_q$ with the bilinear form
$$f(X,Y)=X_1Y_2^T-X_2Y_1^T, \,\, \, \, {\rm  where}\,\,   X=(X_1,X_2), Y=(Y_1,Y_2),  \, X_i, Y_i\in V(n,q).$$   Let $\O =\PG(2n-1, q)$, where   $|\O|=\frac{q^{2n}-1}{q-1}$.
 Let $M=\Sp(2n,q)$ where $n\ge 2$. Consider the action of $M$ on $\O$.
For any $\a\in \O$,  $\D_s=\{\b\in \O\di f(\a, \b)=0\}, $ where $|\D_s|=\frac{q(q^{2n-2}-1)}{q-1}$.  Let $\G=\G (\O, \D_s)$.  Then
 $\{\a, \b\} \in E(\G)$ if and only if  $f(\a, \b)=0$.

Let $\rho=\diag (A, A^\t)$, where $A$ is  the matrix in Eq(\ref{A1}) and $\t$ is the inverse transpose operation.  Then $\rho\in \Sp(2n,q)$ and $\overline {\rho}=\lg \rho\rg /Z\in \PSp(2n,q)$, where $Z=\lg \rho^{\frac{q^n-1}s}\rg$ for $s=(2,q-1)$. Moreover, $\rho(X)=(X_1A, X_2A^\t)$. Set
\begin{equation} \label{s} \begin{array}{llll} &r=\frac{q^n-1}{q-1}, \quad &X_0=(0, 1, 0,\cdots , 0)\in V(n,q),
\quad \b_{1,i}=\lg (X_0A^i,0)\rg, \\
    &\b_{2,i}=\lg (0, X_0(A^i)^{\t})\rg,  &\g_{j,k}=\lg (X_0A^j,X_0(A^{k})^\t)\rg, \quad l=k-j.\end{array}\end{equation}
Then we deal with two cases $n\ge 3$ and $n=2$, separately.

\begin{lem} \label{s1} Let $M=\Sp(2n,q)$ where $n\ge 3$,  acting on the set    $\O$  of (isotropic) points. Then $\G(\O, \Del_s)$  has a H-cycle.
\end{lem}
\demo  Set $s=(2, q-1)$. Using the notations in Eq(\ref{s}), we set
 $$\begin{array}{lll} &B_1=\{ \b_{1,i} \di  0\le i\le r-1\}, \quad   &B_2=\{\b_{2,i}  \di 0\le i\le  r-1\},\\
 &C_l=\{ \g_{j,j+l} \di  0\le j\le \frac{r(q-1)}s-1 \},\,\,  & 0\le |l|\le  sr-1.\end{array}$$
Then $\O$ is a disjoint union of   $\lg \orho\rg $-orbits $B_1, B_2$ and $C_l$, that is
$$\O=B_1\cup B_2\cup  C_0\cup C_1\cup \cdots \cup C_{sr-1}.$$

(1) First, let  $q$ be odd so that $s=2$. Clearly $\G(B_1)$ and  $\G(B_2)$ are complete graphs.
Since $X_0A^{-t}X_0^T-X_0A^{-t}X_0^T=0$ for any $t$,  the induced subgraph  $\G(\{\g_{j-t,k+t}\di \{ 0\le t\le r-1\})$ is a complete graph, for any given $j,k$.
We carry out the proof by the following three steps (i)-(iii).
\vskip 2mm

(i) Let $j=k=0$. Then we have a H-path of $\G(\{\g_{-t, t}\di 0\le t\le r-1\})$:
$$U_0:\quad \g_{-1, 1}, \g_{0,0}, \g_{-2, 2}, \cdots, \g_{-(r-1), (r-1)},$$
whose vertices are respectively contained in  $\lg \orho\rg$-orbits $C_0,C_2, \cdots , C_{2(r-1)}$.
Note that $\g_{-(r-1), (r-1)}\thicksim \g_{0, 2}=\g_{-1,1}^{\orho}$, we get a H-cycle for $\G(\cup_{even\, j} C_j)$:
$$U:\quad U_0,U_0^{\orho},U_0^{\orho^2},\cdots,U_0^{\orho^{\frac {r(q-1)}{2}-1}},\g_{-1,1}.$$

\vskip 2mm
(ii) Let $(j,k)=(0,1)$. Then we have a H-path of $\G(\{\g_{-t, 1+t}\di \{ 0\le t\le r-1\})$:
$$W_0:\quad \g_{0, 1}, \g_{-1,2},\g_{-2,3},\cdots, \g_{-(r-1), r},$$
whose vertices are respectively contained in  $\lg \orho\rg $-orbits $C_1, C_3, \cdots , C_{1+2(r-1)}$.
Since  $\g_{-(r-1), r}\thicksim \g_{1, 2}=\g_{0,1}^{\orho}$, we get
 a H-cycle $W$ for $\G(\cup_{odd\, j} C_j)$:
$$W:\quad W_0,W_0^{\orho},W_0^{\orho^2},\cdots,W_0^{\orho^{\frac {r(q-1)}{2}-1}}, \g_{0,1}.$$

 \vskip 2mm
 (iii)  Choose a H-path $\g_{0,0}\to\g_{-1,1}$   in the H-cycle $U$ (removing the edge $\{\g_{0,0}, \g_{-1,1}\}$);  a H-path $\g_{0,1}\to\g_{-1,2}$  in the H-cycle $W$;
 a H-path $\b_{i,0} \to \b_{i,2}$ in     $\G(B_i)\cong K_r$  where $i=1, 2$.
 Attaching these H-paths,   we obtain a H-cycle $C$  of $\G$ (see Figure 2):
$$C: \quad \g_{0,0}\to\g_{-1,1},\b_{1,0} \to \b_{1,2},\g_{0,1}\to\g_{-1,2},\b_{2,0} \to \b_{2,1},\g_{0,0}.$$

\begin{figure}[htbp]
 \vspace{-2mm}  \centering
  \begin{tikzpicture}[
    scale=0.75, transform shape,
    p/.style={draw, circle, fill=white, inner sep=0pt, minimum size=2.6mm},
    d/.style={fill=white, inner sep=2pt},
    box/.style={draw, rounded corners=6pt},
    lab/.style={anchor=north, yshift=3pt, font=\footnotesize},
    tlab/.style={anchor=south, font=\small},
    midarrow/.style={
      postaction=decorate,
      decoration={markings, mark=at position 0.5 with {\arrow{angle 45}}}
    },
    remidarrow/.style={
      postaction=decorate,
      decoration={markings, mark=at position 0.5 with {\arrowreversed{angle 45}}}
    }
  ]
    % 111
    \node[p] (a10)  at (0, 0)    {};
    \node[d] (d10)  at (2, 0)  {$\cdots$};
    \node[p] (a20)  at (4, 0)  {};
    \node[p] (a30)  at (6, 0)  {};
    \node[p] (c1)   at (8, 0)  {};
    \node[d] (d)    at (8.5, 0)   {$\cdots$};
    \node[p] (c2)   at (9,0)  {};
    \node[p] (b10)  at (11, 0)  {};
    \node[p] (b20)  at (13, 0)  {};
    \node[d] (d20)  at (15, 0)  {$\cdots$};
    \node[p] (b30)  at (17, 0)  {};
    \draw[remidarrow] (a10.east) -- (d10.west);
    \draw[remidarrow] (d10.east) -- (a20.west);
    %\draw[remidarrow] (a20.east) -- (a30.west);
    %\draw[remidarrow] (b10.east) -- (b20.west);
    \draw[remidarrow] (b20.east) -- (d20.west);
    \draw[remidarrow] (d20.east) -- (b30.west);

    % 222
    \node[p] (a11) at (0, -1)  {};
    \node[d] (d11) at (2, -1)  {$\cdots$};
    \node[p] (a21) at (4, -1)  {};
    \node[p] (a31) at (6, -1)  {};
    \node[p] (d1)  at (8, -1.2)  {};
    \node[d] (dd)  at (8.5, -1.2) {$\cdots$};
    \node[p] (d2)  at (9, -1.2)  {};
    \node[p] (b11) at (11, -1)  {};
    \node[p] (b21) at (13, -1)  {};
    \node[d] (d21) at (15, -1)  {$\cdots$};
    \node[p] (b31) at (17, -1)  {};
    \draw[remidarrow] (a11.east) -- (d11.west);
    \draw[remidarrow] (d11.east) -- (a21.west);
    \draw[remidarrow] (a21.east) -- (a31.west);
    \draw[remidarrow] (b11.east) -- (b21.west);
    \draw[remidarrow] (b21.east) -- (d21.west);
    \draw[remidarrow] (d21.east) -- (b31.west);

    % ...
    \node[d] (v1) at (0, -2)   {$\vdots$};
    \node[d] (v2) at (2, -2) {$\vdots$};
    \node[d] (v3) at (4, -2) {$\vdots$};
    \node[d] (v4) at (6, -2)  {$\vdots$};
    \node[d] (v5) at (11, -2)   {$\vdots$};
    \node[d] (v6) at (13, -2) {$\vdots$};
    \node[d] (v7) at (15, -2) {$\vdots$};
    \node[d] (v8) at (17, -2)  {$\vdots$};

    % 333
    \node[p] (a12) at (0, -3)   {};
    \node[d] (d12) at (2, -3) {$\cdots$};
    \node[p] (a22) at (4, -3) {};
    \node[p] (a32) at (6, -3) {};
    \node[p] (b12) at (11, -3)   {};
    \node[p] (b22) at (13, -3) {};
    \node[d] (d22) at (15, -3) {$\cdots$};
    \node[p] (b32) at (17, -3) {};
    \draw[remidarrow] (a12.east) -- (d12.west);
    \draw[remidarrow] (d12.east) -- (a22.west);
    \draw[remidarrow] (a22.east) -- (a32.west);
    \draw[remidarrow] (b12.east) -- (b22.west);
    \draw[remidarrow] (b22.east) -- (d22.west);
    \draw[remidarrow] (d22.east) -- (b32.west);

    % lianjie
    \draw[midarrow] (a10.east)  -- (a31.west);
    \draw[remidarrow] (b30.west)  -- (b11.east);

    \draw[midarrow] (a11.west) -- ($(a11)!0.36!(a32)$);
    \draw[midarrow] ($(a11)!0.64!(a32)$) -- (a32.east);
    \draw[remidarrow] (b31.west) -- ($(b31)!0.36!(b12)$);
    \draw[remidarrow] ($(b31)!0.64!(b12)$) -- (b12.east);

    \draw[midarrow, rounded corners=6pt]
      (a12.west) -- (-1, -3) -- (-1, 0.7) -- (6, 0.7) -- (a30.north);
    \draw[remidarrow, rounded corners=6pt]
      (b32.south) -- (17, -4.1) -- (10, -4.1) -- (10, 0) -- (b10.west);

    \draw[midarrow] (a30.east) -- (c1.west);
    \draw[midarrow, rounded corners=6pt]
      (c2.north) -- (9, 0.7) -- (11, 0.7) -- (b10.north);
    \draw[midarrow] (b20.west) -- (d2.east);
    \draw[midarrow] (d1.west) -- (a20.east);

    % block
    \draw[box] (-0.7, 0.5) rectangle (0.7, -3.9);
    \draw[box] (1.3, 0.5) rectangle (2.7, -3.9);
    \draw[box] (3.3, 0.5) rectangle (4.7, -3.9);
    \draw[box] (5.3, 0.5) rectangle (6.7, -3.9);
    \draw[box] (7.3, 0.5) rectangle (9.7, -0.5);
    \draw[box] (7.3, -0.7) rectangle (9.7, -1.7);
    \draw[box] (10.3, 0.5) rectangle (11.7, -3.9);
    \draw[box] (12.3, 0.5) rectangle (13.7, -3.9);
    \draw[box] (14.3, 0.5) rectangle (15.7, -3.9);
    \draw[box] (16.3, 0.5) rectangle (17.7, -3.9);

    % sign
    \node[lab] at (a10.south)  {$\g_{-(r-1),r-1}$};
    \node[lab] at (a20.south)  {$\g_{0,0}$};
    \node[lab] at (a30.south)  {$\g_{-1,1}$};
    \node[lab] at (b10.south)  {$\g_{0,1}$};
    \node[lab] at (b20.south)  {$\g_{-1.2}$};
    \node[lab] at (b30.south)  {$\g_{-(r-1),r}$};

    \node[lab] at (a11.south) {$\g_{-(r-1),r-1}^{\orho}$};
    \node[lab] at (a21.south) {$\g_{0,0}^{\orho}$};
    \node[lab] at (a31.south) {$\g_{-1,1}^{\orho}$};
    \node[lab] at (b11.south)  {$\g_{0,1}^{\orho}$};
    \node[lab] at (b21.south)  {$\g_{-1.2}^{\orho}$};
    \node[lab] at (b31.south)  {$\g_{-(r-1),r}^{\orho}$};

    \node[lab] at (a12.south) {$\g_{-(r-1),r-1}^{\orho^{\frac{r(q-2)}{2}-1}}$};
    \node[lab] at (a22.south) {$\g_{0,0}^{\orho^{\frac{r(q-2)}{2}-1}}$};
    \node[lab] at (a32.south) {$\g_{-1,1}^{\orho^{\frac{r(q-2)}{2}-1}}$};
    \node[lab] at (b12.south)  {$\g_{0,1}^{\orho^{\frac{r(q-2)}{2}-1}}$};
    \node[lab] at (b22.south)  {$\g_{-1.2}^{\orho^{\frac{r(q-2)}{2}-1}}$};
    \node[lab] at (b32.south)  {$\g_{-(r-1),r}^{\orho^{\frac{r(q-2)}{2}-1}}$};

    \node[lab] at (c1.south) {$\b_{1,0}$};
    \node[lab] at (c2.south) {$\b_{1.2}$};
    \node[lab] at (d1.south) {$\b_{2,1}$};
    \node[lab] at (d2.south) {$\b_{2,0}$};

    \node[tlab] at (0, 0.6) {$C_{2(r-1)}$};
    \node[tlab] at (2, 0.6) {$\cdots$};
    \node[tlab] at (4, 0.6) {$C_{0}$};
    \node[tlab] at (6, 0.6) {$C_{2}$};
    \node[tlab] at (11, 0.6) {$C_{1}$};
    \node[tlab] at (13, 0.6) {$C_{3}$};
    \node[tlab] at (15, 0.6) {$\cdots$};
    \node[tlab] at (17, 0.6) {$C_{2r-1}$};
    \node[tlab] at (8, 0.6) {$B_1$};
    \node[tlab] at (8, -2.3) {$B_2$};
  \end{tikzpicture}
  \caption{{\footnotesize  H-cycle $C$ of $\G$ in Lemma 5.1.(1).}}
\end{figure}

\vskip 2mm

(2)  Second, let  $q$ be even so that $s=1$.
 Then $\G(\{\g_{j-\frac t{2},k+\frac t{2}}\di 0\le t\le r-1\})$ is a complete graph.
 Let $j=k=0$. Then we have a H-path of $\G(\{\g_{-\frac t2, \frac t2}\di 0\le t\le r-1\})$:

$$U_0:\quad \g_{-\frac 1{2},\frac 1{2}},\g_{0,0},\g_{-1,1},\cdots,\g_{-\frac {r-1}{2},\frac {r-1}{2}},$$
whose vertices are respectively contained in  $\lg \rho \rg$-orbits $C_0,C_1,C_2,\cdots,C_{r-1}$.
Check $\g_{-\frac {r-1}{2},\frac {r-1}{2}} \thicksim \g_{\frac 1{2},\frac 3{2}} = \g_{-\frac 1{2},\frac 1{2}}^{\rho}$. Then we get a H-cycle $U$ for $\G(\cup_{j}C_j)$:
$$U:\quad U_0,U_0^\rho,U_0^{\rho^2},\cdots,U_0^{\rho^{r(q-1)-1}},\g_{-\frac 1{2},\frac 1{2}}.$$
Pick up a H-path $\g_{-\frac 1{2},\frac 1{2}}\to\g_{0,0}$ in $U$, a H-path $\b_{1,1} \to \b_{1,0}$ in  $\G(B_1)$  and a H-path  $\b_{2,1} \to \b_{2,\frac 1{2}}$   in $\G(B_2)$.
Attaching these H-paths, we obtain a H-cycle $C$ of $\G$:
$$C: \quad \g_{-\frac 1{2},\frac 1{2}}\to\g_{0,0},\b_{1,1} \to \b_{1,0},\b_{2,1} \to \b_{2,\frac 1{2}},\g_{-\frac 1{2},\frac 1{2}}.$$
\qed

\begin{lem} \label{s2} Let $M=\Sp(4,q)$,  acting on the set    $\O$  of (isotropic) points. Then $\G(\O, \Del_s)$  has a H-cycle.
\end{lem}
\demo  Use the notations in Eq(\ref{s}). Now,  $$A=\left(\begin{array}{cc} 0 &1\\a_0 &a_1\end{array}\right),\, \,  A^{-1}=\left(\begin{array}{cc} -a_1a_0^{-1} &a_0^{-1}\\1 &0\end{array}\right), s=(2, q-1),  \, \, {\rm and} \, \, X_0=(0,1).$$
 Set
$$\begin{array}{llll}
B_1&=&\{\b_{1,i}\di 0\le i \le r-1\}, \quad  &\b_{1,i}=\lg(X_0A^i,0_m)\rg ;\\
B_2&=&\{\b_{2,i}\di 0\le i \le r-1\}, \quad  &\b_{2,i}=\lg(0_m,X_0(A^i)^\t)\rg ;\\
C_1&=&\{\d_{1, k}\di 0\le k \le \frac {q^{n}-1}{s}-1\}, &\d_{1, k}=\lg(X_0A^k,(1,0)(A^k)^\t)\rg ; \\
C_{\th}&=&\{\d_{\th, k}\di 0\le k \le \frac {q^{n}-1}{2}-1\}, & \d_{\th, k}=\lg(X_0A^k,(\th,0)(A^k)^\t)\rg , \, {only\, if}\,  s=2; \\
D_l&=&\{\g_{l, k}\di 0\le k \le \frac {q^n-1}{s}-1\}, &\g_{l, k}=\lg(X_0A^k,(l,1)(A^k)^\t)\rg,\,\,l\in \FF_q;\\
E_l&=&\{\ph_{l, k}\di 0\le k \le \frac {q^n-1}{2}-1\}, &\ph_{l, k}=\lg(X_0A^k,(l,\th)(A^k)^\t)\rg,\,\,l\in \FF_q, \, {only\, if}\\
 &&& \hskip 5.2cm s=2.
\end{array} $$
Check that   $\O$ is a disjoint union of  the above $\lg\rho\rg$-orbits, that is
$$\begin{array}{ll} \O=B_1\cup B_2\cup  C_1\cup  C_\th \cup (\cup _{l\in \FF_q} D_l)\cup (\cup_{l\in \FF_q}E_l), \, {\rm for }\, s=2;\\
\O=B_1\cup B_2\cup  C_1\cup (\cup _{l\in \FF_q} D_l),  \, {\rm for }\, s=1.
   \end{array}$$
Moreover,   $\G(B_1)$, $\G(B_2)$, $\G(\{\g_{l,k}\di l\in \FF_q\})$ and $\G(\{\ph_{l,k}\di l\in \FF_q\})$ for any given  $k$, are complete graphs.
\vskip 2mm
(1) First, let $s=2$. Recall $a_0$ and  $a_1$ are two entries of $A$.
There is at least one of  $\{a_0, \th a_0 \}$ which is not $a_1$, no loss, we assume $\th a_0 \ne a_1$. Carry out the proof by three steps below.
\vskip 2mm
(i)  Noting   $\G(\{\g_{l,1}\di l\in \FF_q\})$ is complete,
we may pick up a H-path   $U_0$ of  $\G(\{\g_{l,1}\di l\in \FF_q\}\cup\{\d_{1,0},\d_{\th,0}\})$:
$$\begin{array} {lll} U_0&=&\g_{a_1,1}, \g_{a_0,1}, \d_{1,0}, \d_{\th,0}, \g_{\th a_0,1}, \cdots , \g_{c,1}, \g_{0,1},\quad  {\rm if}\,  a_1\ne a_0; {\rm or}\\
 U_0&=&\g_{a_1,1}, \d_{1,0}, \d_{\th,0}, \g_{\th a_0,1}, \cdots , \g_{c,1}, \g_{0,1},\quad  {\rm if}\,  a_1=a_0,
\end{array} $$
for some $c\ne 0, a_0, a_1$. Note that $a_1\ne 0$. Since $\g_{0,1}\thicksim\g_{a_1,2}=\g_{a_1,1}^{\orho}$, we get a H-path $U$ in $\G((\cup_lD_l)\cup C_1\cup C_\th)$ (see Figure 3):
$$U:\quad \g_{0,1},U_0^{\orho},U_0^{\orho^2},\cdots,U_0^{\orho^{\frac{q^2-1}{2}-1}},\g_{a_1,1},\cdots,\g_{c,1}.$$
\vskip 2mm
(ii) In the complete graph $\G(\{\ph_{l,0}\di l\})$, pick a H-path:
$$W_0:\quad \ph_{\th a_1,0},\cdots,\ph_{d,0},\ph_{0,0},$$
 for some  $d\ne 0, \th a_1.$   Since $\ph_{0,0}\thicksim\ph_{\th a_1,1}=\ph_{\th a_1,0}^{\orho}$, we get a H-path $W$ in $\G(\cup_l E_l)$ (see the left part of Figure 4):
$$W:\quad \ph_{0,0},W_0^{\orho},W_0^{\orho^2},\cdots,W_0^{\orho^{\frac{q^2-1}{2}-1}},\ph_{\th a_1,0},\cdots,\ph_{d,0}.$$

\vskip 2mm
(iii) Taking  a H-path $\b_{2,2}\to\b_{2,1}$ in $\G(B_2)$ and a H-path $\b_{1, r-1}\to\b_{1,0}$ in $\G(B_1)$, we finally get a H-cycle $C$ of $\G$ (see Figure 4):
$$C:  \quad U,\b_{2,2}\to\b_{2,1},W^{-1},\b_{1,r-1}\to\b_{1,0}, \g_{0,1}.$$

\begin{figure}[htbp]
  \vspace{-2mm}
  \centering
  \begin{tikzpicture}[
    scale=0.75, transform shape,
    p/.style={draw, circle, fill=white, inner sep=0pt, minimum size=2.6mm},
    d/.style={fill=white, inner sep=2pt},
    box/.style={draw, rounded corners=6pt},
    lab/.style={anchor=north, yshift=3pt, font=\footnotesize},
    tlab/.style={anchor=south, font=\small},
    midarrow/.style={
      postaction=decorate,
      decoration={markings, mark=at position 0.5 with {\arrow{angle 45}}}
    },
    remidarrow/.style={
      postaction=decorate,
      decoration={markings, mark=at position 0.5 with {\arrowreversed{angle 45}}}
    }
  ]
    % 111
    \node[p] (b10)  at (0, 0)  {};
    \node[p] (a20)  at (2, 0)  {};
    \node[p] (a30)  at (4, 0)  {};
    \node[p] (a40)  at (6, 0)  {};
    \node[p] (a50)  at (8, 0)  {};
    \node[d] (d10)  at (10,0)  {$\cdots$};
    \node[p] (a60)  at (12,0)  {};
    \node[p] (a70)  at (14,0)  {};
    \draw[midarrow] (a10.east) -- (a20.west);
    \draw[midarrow] (a20.east) -- (a30.west);
    \draw[midarrow] (a30.east) -- (a40.west);
    \draw[midarrow] (a40.east) -- (a50.west);
    \draw[midarrow] (a50.east) -- (d10.west);
    \draw[midarrow] (d10.east) -- (a60.west);
    %\draw[midarrow] (a60.east) -- (a70.west);
    % 222
    \node[p] (a11)  at (0, -1)  {};
    \node[p] (a21)  at (2, -1)  {};
    \node[p] (a31)  at (4, -1)  {};
    \node[p] (a41)  at (6, -1)  {};
    \node[p] (a51)  at (8, -1)  {};
    \node[d] (d11)  at (10,-1)  {$\cdots$};
    \node[p] (a61)  at (12,-1)  {};
    \node[p] (a71)  at (14,-1)  {};
    \draw[midarrow] (a11.east) -- (a21.west);
    \draw[midarrow] (a21.east) -- (a31.west);
    \draw[midarrow] (a31.east) -- (a41.west);
    \draw[midarrow] (a41.east) -- (a51.west);
    \draw[midarrow] (a51.east) -- (d11.west);
    \draw[midarrow] (d11.east) -- (a61.west);
    \draw[midarrow] (a61.east) -- (a71.west);
    % ...
    \node[d] (v1)  at (0, -2)  {$\vdots$};
    \node[d] (v2)  at (2, -2)  {$\vdots$};
    \node[d] (v3)  at (4, -2)  {$\vdots$};
    \node[d] (v4)  at (6, -2)  {$\vdots$};
    \node[d] (v5)  at (8, -2)  {$\vdots$};
    \node[d] (v6)  at (10,-2)  {$\vdots$};
    \node[d] (v7)  at (12,-2)  {$\vdots$};
    \node[d] (v8)  at (14,-2)  {$\vdots$};
    % 333
    \node[p] (a12)  at (0, -3)  {};
    \node[p] (a22)  at (2, -3)  {};
    \node[p] (a32)  at (4, -3)  {};
    \node[p] (a42)  at (6, -3)  {};
    \node[p] (a52)  at (8, -3)  {};
    \node[d] (d12)  at (10,-3)  {$\cdots$};
    \node[p] (a62)  at (12,-3)  {};
    \node[p] (a72)  at (14,-3)  {};
    \draw[midarrow] (a12.east) -- (a22.west);
    \draw[midarrow] (a22.east) -- (a32.west);
    \draw[midarrow] (a32.east) -- (a42.west);
    \draw[midarrow] (a42.east) -- (a52.west);
    \draw[midarrow] (a52.east) -- (d12.west);
    \draw[midarrow] (d12.east) -- (a62.west);
    \draw[midarrow] (a62.east) -- (a72.west);

    % lian
    \draw[midarrow] (a70.west) -- (a11.east);
    \draw[midarrow] (a71.west) -- ($(a71)!0.36!(a12)$);
    \draw[midarrow] ($(a71)!0.64!(a12)$) -- (a12.east);
    \draw[midarrow, rounded corners=6pt]
    (a72.south) -- (14,-4) -- (-1,-4) -- (-1,0) -- (a10.west);

    % block
    \draw[box] (-0.7, 0.5) rectangle (0.7, -3.8);
    \draw[box] (1.3, 0.5) rectangle (2.7, -3.8);
    \draw[box] (3.3, 0.5) rectangle (4.7, -3.8);
    \draw[box] (5.3, 0.5) rectangle (6.7, -3.8);
    \draw[box] (7.3, 0.5) rectangle (8.7, -3.8);
    \draw[box] (9.3, 0.5) rectangle (10.7, -3.8);
    \draw[box] (11.3, 0.5) rectangle (12.7, -3.8);
    \draw[box] (13.3, 0.5) rectangle (14.7, -3.8);

    % sign
    \node[lab] at (a10.south)  {$\g_{a_1,1}$};
    \node[lab] at (a20.south)  {$\g_{a_0,1}$};
    \node[lab] at (a30.south)  {$\d_{1,0}$};
    \node[lab] at (a40.south)  {$\d_{\th,0}$};
    \node[lab] at (a50.south)  {$\g_{a_0\th,1}$};
    \node[lab] at (a60.south)  {$\g_{c,1}$};
    \node[lab] at (a70.south)  {$\g_{0,1}$};

    \node[lab] at (a11.south)  {$\g_{a_1,2}$};
    \node[lab] at (a21.south)  {$\g_{a_0,2}$};
    \node[lab] at (a31.south)  {$\d_{1,1}$};
    \node[lab] at (a41.south)  {$\d_{\th,1}$};
    \node[lab] at (a51.south)  {$\g_{a_0\th,2}$};
    \node[lab] at (a61.south)  {$\g_{c,2}$};
    \node[lab] at (a71.south)  {$\g_{0,2}$};

    \node[lab] at (a12.south)  {$\g_{a_1,0}$};
    \node[lab] at (a22.south)  {$\g_{a_0,0}$};
    \node[lab] at (a32.south)  {$\d_{1,\frac{q^2-3}{2}}$};
    \node[lab] at (a42.south)  {$\d_{\th,\frac{q^2-3}{2}}$};
    \node[lab] at (a52.south)  {$\g_{a_0\th,0}$};
    \node[lab] at (a62.south)  {$\g_{c,0}$};
    \node[lab] at (a72.south)  {$\g_{0,0}$};

    \node[tlab] at (0, 0.5) {$D_{a_1}$};
    \node[tlab] at (2, 0.5) {$D_{a_0}$};
    \node[tlab] at (4, 0.5) {$C_1$};
    \node[tlab] at (6, 0.5) {$C_{\th}$};
    \node[tlab] at (8, 0.5) {$D_{a_0\th}$};
    \node[tlab] at (10, 0.5) {$\cdots$};
    \node[tlab] at (12, 0.5) {$D_{c}$};
    \node[tlab] at (14, 0.5) {$D_{0}$};
  \end{tikzpicture}
   \caption{{\footnotesize  H-path $U$ of  $\G((\cup_lD_l)\cup C_1\cup C_\th)$ in Lemma 5.2(1).}}

  \vskip 2mm
  \begin{tikzpicture}[
    scale=0.75, transform shape,
    p/.style={draw, circle, fill=white, inner sep=0pt, minimum size=2.6mm},
    d/.style={fill=white, inner sep=2pt},
    box/.style={draw, rounded corners=6pt},
    lab/.style={anchor=north, yshift=3pt, font=\footnotesize},
    tlab/.style={anchor=south, font=\small},
    midarrow/.style={
      postaction=decorate,
      decoration={markings, mark=at position 0.5 with {\arrow{angle 45}}}
    },
    remidarrow/.style={
      postaction=decorate,
      decoration={markings, mark=at position 0.5 with {\arrowreversed{angle 45}}}
    }
  ]
    % 111
    \node[p] (b10)  at (0, 0)  {};
    \node[d] (d10)  at (2,0)  {$\cdots$};
    \node[p] (b20)  at (4, 0)  {};
    \node[p] (b30)  at (6,0)  {};
    \node[p] (b40)  at (8,0)  {};
    \node[d] (d20)  at (9,0)  {$\cdots$};
    \node[p] (b50)  at (10,0)  {};

    \draw[remidarrow] (b10.east) -- (d10.west);
    \draw[remidarrow] (d10.east) -- (b20.west);
    %\draw[midarrow] (b20.east) -- (b30.west);
    % 222
    \node[p] (b11)  at (0, -1)  {};
    \node[d] (d11)  at (2,-1)  {$\cdots$};
    \node[p] (b21)  at (4, -1)  {};
    \node[p] (b31)  at (6,-1)  {};
    \node[p] (b41)  at (8,-1.2)  {};
    \node[d] (d21)  at (9,-1.2)  {$\cdots$};
    \node[p] (b51)  at (10,-1.2)  {};
    \draw[remidarrow] (b11.east) -- (d11.west);
    \draw[remidarrow] (d11.east) -- (b21.west);
    \draw[remidarrow] (b21.east) -- (b31.west);
    % ...
    \node[d] (v1)  at (0, -2)  {$\vdots$};
    \node[d] (v2)  at (2, -2)  {$\vdots$};
    \node[d] (v3)  at (4, -2)  {$\vdots$};
    \node[d] (v4)  at (6, -2)  {$\vdots$};
    % 333
    \node[p] (b12)  at (0, -3)  {};
    \node[d] (d12)  at (2,-3)  {$\cdots$};
    \node[p] (b22)  at (4, -3)  {};
    \node[p] (b32)  at (6,-3)  {};
    \node[p] (b42)  at (8,-2.4)  {};
    \node[d] (d22)  at (9,-2.4)  {$\cdots$};
    \node[p] (b52)  at (10,-2.4)  {};
    \draw[remidarrow] (b12.east) -- (d12.west);
    \draw[remidarrow] (d12.east) -- (b22.west);
    \draw[remidarrow] (b22.east) -- (b32.west);

    % lian
    \draw[remidarrow] (b30.west) -- (b11.east);
    \draw[remidarrow] (b31.west) -- ($(b31)!0.36!(b12)$);
    \draw[remidarrow] ($(b31)!0.64!(b12)$) -- (b12.east);
    \draw[remidarrow, rounded corners=6pt]
    (b32.south) -- (6,-4) -- (-1,-4) -- (-1,0) -- (b10.west);
    \draw[midarrow] (b52.west) -- ($(b52)!0.3!(b42)$);
    \draw[midarrow] ($(b52)!0.7!(b42)$) -- (b42.east);
    \draw[midarrow] (b42.east) -- (b51.west);
    \draw[midarrow] (b51.west) -- ($(b51)!0.3!(b41)$);
    \draw[midarrow] ($(b51)!0.7!(b41)$) -- (b41.east);
    \draw[midarrow] (b41.west) -- (b20.east);
    \draw[midarrow] (b30.east) -- (b40.west);
    \draw[midarrow] (b40.east) -- ($(b40)!0.3!(b50)$);
    \draw[midarrow] ($(b40)!0.7!(b50)$) -- (b50.west);
    \draw[midarrow, rounded corners=6pt]
    (b50.east) -- (10.9, 0) -- (10.9, -2.4) -- (b52.east);

    % block
    \draw[box] (-0.7, 0.5) rectangle (0.7, -3.8);
    \draw[box] (1.3, 0.5) rectangle (2.7, -3.8);
    \draw[box] (3.3, 0.5) rectangle (4.7, -3.8);
    \draw[box] (5.3, 0.5) rectangle (6.7, -3.8);
    \draw[box] (7.3, 0.5) rectangle (10.7,-0.5);
    \draw[box] (7.3, -0.7) rectangle (10.7,-1.7);
    \draw[box] (7.3, -1.9) rectangle (10.7,-2.9);

    % sign
    \node[lab] at (b10.south)  {$\ph_{\th a_1,0}$};
    \node[lab] at (b20.south)  {$\ph_{d,0}$};
    \node[lab] at (b30.south)  {$\ph_{0,0}$};

    \node[lab] at (b11.south)  {$\ph_{\th a_1,1}$};
    \node[lab] at (b21.south)  {$\ph_{d,1}$};
    \node[lab] at (b31.south)  {$\ph_{0,1}$};

    \node[lab] at (b12.south)  {$\ph_{\th a_1,\frac{q^2-3}{2}}$};
    \node[lab] at (b22.south)  {$\ph_{d,\frac{q^2-3}{2}}$};
    \node[lab] at (b32.south)  {$\ph_{0,\frac{q^2-3}{2}}$};

    \node[lab] at (b40.south)  {$\b_{1,r-1}$};
    \node[lab] at (b50.south)  {$\b_{1,0}$};
    \node[lab] at (b41.south)  {$\b_{2,1}$};
    \node[lab] at (b51.south)  {$\b_{2,2}$};
    \node[lab] at (b42.south)  {$\g_{c,1}$};
    \node[lab] at (b52.south)  {$\g_{0,1}$};

    \node[tlab] at (0, 0.5) {$E_{\th a_1}$};
    \node[tlab] at (2, 0.5) {$\cdots$};
    \node[tlab] at (4, 0.5) {$E_{d}$};
    \node[tlab] at (6, 0.5) {$E_{0}$};
    \node[tlab] at (11.2, -0.3) {$B_1$};
    \node[tlab] at (11.2, -1.5) {$B_2$};
    \node[tlab] at (12.4, -2.7) {$(\cup_l D_l)\cup C_1\cup C_\th$};
\end{tikzpicture}
  \caption{{\footnotesize  H-cycle $C$ of $\G$  in Lemma 5.2(1).}}
\end{figure}

\vskip 2mm (2) Second,  let  $s=1$. Then $C_\th$ and all  $E_l$ are disappeared.
For $q=2$, $|\D_s|\ge \frac{1}{3}|\O|$ and so let $q>2$. Since $\d_{1,ir}=\lg(0,\th^i,\th^{-i},0)\rg$,  the graph $\G(\{\d_{1,ir}\di 0\le i \le q-2\})$ is a complete graph. Similarly,   we have  $q+1$ such complete graphs $\G(\{\d_{1,j+ir}\di 0\le i \le q-2\})$, where $0\le j \le q$. Let $\d_{1,i'r}=\lg(0,a'^{-1},a',0)\rg$ for some $i'$, where $a'^2=a_0\ne 1$. Check  $\d_{1,0}\thicksim \d_{1,i'r+1}$. In the complete graph $\G(\{\d_{1,ir}\di 0\le i\le q-2 \})$, pick a path $U_0$ leading $\d_{1,i'r}$ to $\d_{1,0}$, and then get a H-path in $\G(C_1)$:
$$U:\quad U_0,U_0^\rho,U_0^{\rho^2},\cdots,U_0^{\rho^q}.$$
Note that $a_1 \ne 0$.
In the complete graph $\G(\{\g_{l,0}\di l\in \FF_q\})$, pick a path $W_0$: $\g_{a_1,0},\cdots,\g_{b,0},\g_{0,0}$, $b\ne 0, a_1.$   Since $\g_{0,0}\thicksim\g_{a_1,1}$, we get a H-path in $\G(\cup_lD_l)$:
$$W:\quad \g_{0,0},W_0^\rho,W_0^{\rho^2},\cdots,W_0^{\rho^{r(q-1)-1}},\g_{a_1,0},\cdots,\g_{b,0}.$$
Take a H-path $\b_{2,0}\to\b_{2,1}$ in $\G(B_2)$ and a H-path $\b_{1,r-1}\to\b_{1,0}$ in $\G(B_1)$. Then we finally get a H-cycle $C$ of  $\G$ (see Figure 5.):
$$C:\quad U,\b_{2,0}\to\b_{2,1},W^{-1},\b_{1,r-1}\to\b_{1,0},\d_{0,i'r}.$$
\vspace{-2mm}  \qed
\begin{figure}[htbp]
 \vspace{-2mm} \centering
  \begin{tikzpicture}[
    scale=0.75, transform shape,
    p/.style={draw, circle, fill=white, inner sep=0pt, minimum size=2.6mm},
    d/.style={fill=white, inner sep=2pt},
    box/.style={draw, rounded corners=6pt},
    lab/.style={anchor=north, yshift=3pt, font=\footnotesize},
    tlab/.style={anchor=south, font=\small},
    midarrow/.style={
      postaction=decorate,
      decoration={markings, mark=at position 0.5 with {\arrow{angle 45}}}
    },
    remidarrow/.style={
      postaction=decorate,
      decoration={markings, mark=at position 0.5 with {\arrowreversed{angle 45}}}
    }
  ]
    % 111
    \node[p] (b10)  at (0, 0)  {};
    \node[d] (d10)  at (2,0)  {$\cdots$};
    \node[p] (b20)  at (4, 0)  {};
    \node[p] (b30)  at (6,0)  {};
    \node[p] (b40)  at (8,0)  {};
    \node[d] (d20)  at (9,0)  {$\cdots$};
    \node[p] (b50)  at (10,0)  {};
    \node[p] (b60)  at (12,0)  {};
    \node[d] (d30)  at (13,0)  {$\cdots$};
    \node[p] (b70)  at (14,0)  {};

    \draw[remidarrow] (b10.east) -- (d10.west);
    \draw[remidarrow] (d10.east) -- (b20.west);
    %\draw[midarrow] (b20.east) -- (b30.west);
    % 222
    \node[p] (b11)  at (0, -1)  {};
    \node[d] (d11)  at (2,-1)  {$\cdots$};
    \node[p] (b21)  at (4, -1)  {};
    \node[p] (b31)  at (6,-1)  {};
    \node[p] (b41)  at (8,-1.2)  {};
    \node[d] (d21)  at (9,-1.2)  {$\cdots$};
    \node[p] (b51)  at (10,-1.2)  {};
    \node[p] (b61)  at (12,-1)  {};
    \node[d] (d31)  at (13,-1)  {$\cdots$};
    \node[p] (b71)  at (14,-1)  {};
    \draw[remidarrow] (b11.east) -- (d11.west);
    \draw[remidarrow] (d11.east) -- (b21.west);
    \draw[remidarrow] (b21.east) -- (b31.west);
    % ...
    \node[d] (v1)  at (0, -2)  {$\vdots$};
    \node[d] (v2)  at (2, -2)  {$\vdots$};
    \node[d] (v3)  at (4, -2)  {$\vdots$};
    \node[d] (v4)  at (6, -2)  {$\vdots$};
    \node[d] (v5)  at (12, -2)  {$\vdots$};
    \node[d] (v6)  at (14, -2)  {$\vdots$};
    % 333
    \node[p] (b12)  at (0, -3)  {};
    \node[d] (d12)  at (2,-3)  {$\cdots$};
    \node[p] (b22)  at (4, -3)  {};
    \node[p] (b32)  at (6,-3)  {};
    \node[p] (b62)  at (12,-3)  {};
    \node[d] (d32)  at (13,-3)  {$\cdots$};
    \node[p] (b72)  at (14,-3)  {};
    \draw[remidarrow] (b12.east) -- (d12.west);
    \draw[remidarrow] (d12.east) -- (b22.west);
    \draw[remidarrow] (b22.east) -- (b32.west);

    % lian
    \draw[remidarrow] (b30.west) -- (b11.east);
    \draw[remidarrow] (b31.west) -- ($(b31)!0.36!(b12)$);
    \draw[remidarrow] ($(b31)!0.64!(b12)$) -- (b12.east);
    \draw[remidarrow, rounded corners=6pt]
    (b32.south) -- (6,-4) -- (-1,-4) -- (-1,0) -- (b10.west);
    \draw[midarrow] (b60.east) -- ($(b60)!0.3!(b70)$);
    \draw[midarrow] ($(b60)!0.7!(b70)$) -- (b70.west);
    \draw[midarrow] (b70.west) -- (b61.east);
    \draw[midarrow] (b61.east) -- ($(b61)!0.3!(b71)$);
    \draw[midarrow] ($(b61)!0.7!(b71)$) -- (b71.west);
    \draw[midarrow] (b71.west) -- ($(b71)!0.36!(b62)$);
    \draw[midarrow] ($(b71)!0.64!(b62)$) -- (b62.east);
    \draw[midarrow] (b62.east) -- ($(b62)!0.3!(b72)$);
    \draw[midarrow] ($(b62)!0.7!(b72)$) -- (b72.west);
    \draw[midarrow, rounded corners=6pt]
    (b72.south) -- (14,-4) -- (10,-4) -- (b51.south);
    \draw[midarrow] (b51.west) -- ($(b51)!0.3!(b41)$);
    \draw[midarrow] ($(b51)!0.7!(b41)$) -- (b41.east);
    \draw[midarrow] (b41.west) -- (b20.east);
    \draw[midarrow] (b30.east) -- (b40.west);
    \draw[midarrow] (b40.east) -- ($(b40)!0.3!(b50)$);
    \draw[midarrow] ($(b40)!0.7!(b50)$) -- (b50.west);
    \draw[midarrow] (b50.east) -- (b60.west);

    % block
    \draw[box] (-0.7, 0.5) rectangle (0.7, -3.8);
    \draw[box] (1.3, 0.5) rectangle (2.7, -3.8);
    \draw[box] (3.3, 0.5) rectangle (4.7, -3.8);
    \draw[box] (5.3, 0.5) rectangle (6.7, -3.8);
    \draw[box] (7.3, 0.5) rectangle (10.7,-0.5);
    \draw[box] (7.3, -0.7) rectangle (10.7,-1.7);
    \draw[box] (11.3, 0.5) rectangle (14.7,-3.8);

    % sign
    \node[lab] at (b10.south)  {$\g_{a_1,0}$};
    \node[lab] at (b20.south)  {$\g_{b,0}$};
    \node[lab] at (b30.south)  {$\g_{0,0}$};

    \node[lab] at (b11.south)  {$\g_{a_1,1}$};
    \node[lab] at (b21.south)  {$\g_{b,1}$};
    \node[lab] at (b31.south)  {$\g_{0,1}$};

    \node[lab] at (b12.south)  {$\g_{\th a_1,q^2-2}$};
    \node[lab] at (b22.south)  {$\g_{b,q^2-2}$};
    \node[lab] at (b32.south)  {$\g_{0,q^2-2}$};

    \node[lab] at (b40.south)  {$\b_{1,r-1}$};
    \node[lab] at (b50.south)  {$\b_{1,0}$};
    \node[lab] at (b41.south)  {$\b_{2,1}$};
    \node[lab] at (b51.south)  {$\b_{2,0}$};

    \node[lab] at (b60.south)  {$\d_{1,i'r}$};
    \node[lab] at (b70.south)  {$\d_{1,0}$};
    \node[lab] at (b61.south)  {$\d_{1,i'r+1}$};
    \node[lab] at (b71.south)  {$\d_{1,1}$};
    \node[lab] at (b62.south)  {$\d_{1,i'r+q}$};
    \node[lab] at (b72.south)  {$\d_{1,q}$};

    \node[tlab] at (0, 0.5) {$D_{a_1}$};
    \node[tlab] at (2, 0.5) {$\cdots$};
    \node[tlab] at (4, 0.5) {$D_{b}$};
    \node[tlab] at (6, 0.5) {$D_{0}$};
    \node[tlab] at (8, 0.5) {$B_1$};
    \node[tlab] at (8, -2.3) {$B_2$};
    \node[tlab] at (12, 0.5) {$C_1$};

  \end{tikzpicture}
  \caption{{\footnotesize H-cycle $C$ of $\G$ in Lemma 5.2.(2).}}
\vspace{-2mm}
\end{figure}

\subsection{Action  on  maximal totally isotropic subspaces}

\begin{lem}\label{s2}
Let $M=\Sp(4,q)$, acting on   the set $\O$ of  maximal totally isotropic subspaces. Then  $\G(\O, \Del_s)$  has a H-cycle.
\end{lem}
\demo Use the symplectic form  $f$  given in last subsection.
Set $V=\lg u_1,u_2\rg\perp\lg v_1,v_2\rg$, where $\lg u_1,u_2\rg$ and $\lg v_1,v_2\rg$ are hyperbolic planes, that is $f(u_1,u_2)=f(v_1, v_2)=1$.
Then $|\O|=q^3+q^2+q+1$.
For any $U\in \O$,
$\D_s=\{W\in\O \di \dim(W\cap U)=1\},$  where $|\D_s|=q(q+1).$  Let $\G=\G(\O,\D_s)$.  Then
 $\{W_1, W_2\} \in E(\G)$ if and only if  $\dim(W_1\cap W_2)=1$.

For any $v\in V$, set $\o v=\{W\in \O \di \lg v \rg \in W\}$ so that   $|\o v|=q+1$. Let $\{v_i|0\le i\le q\}$ be all points in $\lg v_1,v_2\rg$.
Let $\FF_q^*=\lg \th \rg$. Pick up an  $u\in\lg u_1,u_2\rg$, and set
$\d_{i0}=v_i$ and  $\d_{ij}=u+\th^j v_i$, where $0\le i \le q$ and  $1\le j \le q-1$. Then  $\O$ is a union of
all $\o u$ and  all $\o \d_{ij}$, that is
$$\O=(\cup_{u\in \lg u_1, u_2\rg } \ou ) \cup (\cup_{0\le i \le q,  0 \le j \le q-1}  \o \d_{ij}).$$
\f Check that for any given $i$, $\o \d_{ij} \cap \o \d_{ij'}=\lg u,v_i\rg$; and for any given $i_1 \neq i_2$, there is a matching between $\o \d_{i_1,j_1}$ and $\o\d_{i_2,j_2}$ such that
$$\lg \d_{i_1,j_1},w \rg\to \lg \d_{i_2,j_2}, \,\,  w \rg, \forall{w} \in \d_{i_1,j_1}^\perp\cap \d_{i_2,j_2}^\perp.$$
The the induced subgraphs $X(\o u)$ and $X(\o \d_{ij})$ are complete.   Take a sequence $\mathcal L$  of $q(q+1)$ of subgraphs $X(\o \d_{ij})$:
$$\mathcal L: \o \d_{00},\o\d_{10},\cdots,\o\d_{q0};\o\d_{01},\o\d_{11},\cdots,\o\d_{q1};\cdots;\o\d_{0,q-1},\o\d_{1,q-1},\cdots,\o\d_{q,q-1}.$$
Along  the sequence $\mathcal L$, choose a H-path $\mathcal H_{ij}$ in $\o\d_{ij}\setminus \lg u,v_i\rg$, such that the ending vertex of $\mathcal H_{ij}$ is adjacent to the starting vertex of the subsequent  term, by using the above matching. Let $\mathcal H_{\infty}$ be a H-path in $\o u$ leading  $\lg u,v_q\rg$ to $\lg u,v_0\rg$. Then we get a H-cycle $C$ of $\G$ (see Figure 6):
$$C: \quad \lg u,v_0\rg,\mathcal H_{00},\mathcal H_{10},\cdots,\mathcal H_{q0};\cdots;\mathcal H_{0,q-1},\mathcal H_{1,q-1},\cdots,\mathcal H_{q,q-1};\mathcal H_{\infty}.$$\qed

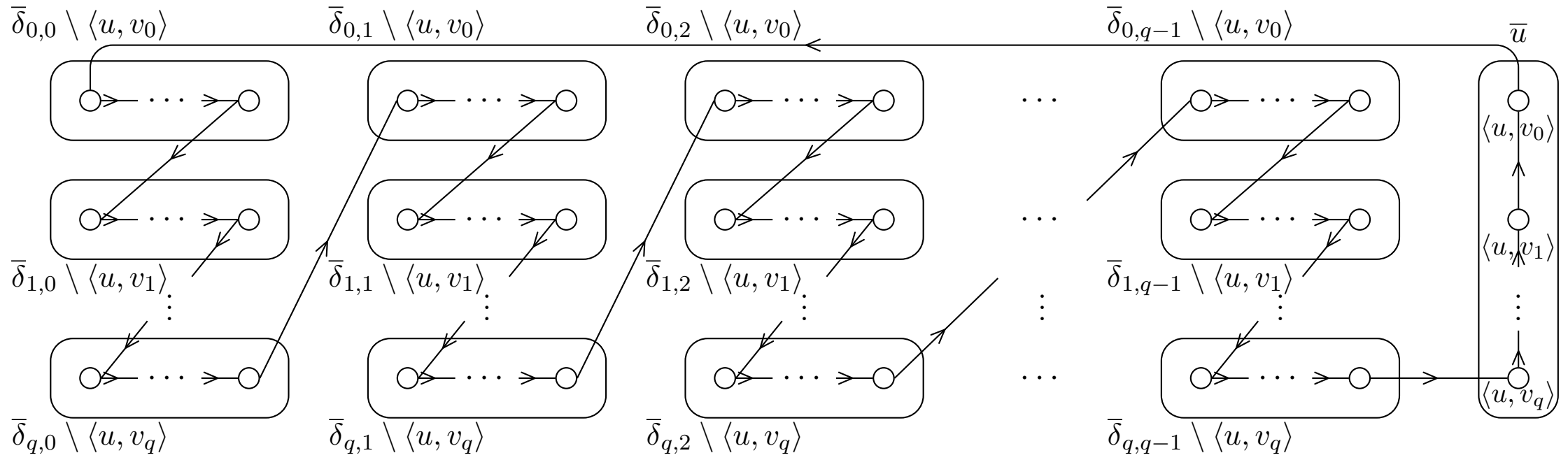
\begin{figure}[htbp]
 \vspace{-2mm} \centering
  \begin{tikzpicture}[
    scale=0.75, transform shape,
    p/.style={draw, circle, fill=white, inner sep=0pt, minimum size=2.6mm},
    d/.style={fill=white, inner sep=2pt},
    box/.style={draw, rounded corners=6pt},
    lab/.style={anchor=north, yshift=3pt, font=\footnotesize},
    tlab/.style={anchor=south, font=\small},
    midarrow/.style={
      postaction=decorate,
      decoration={markings, mark=at position 0.5 with {\arrow{angle 45}}}
    },
    remidarrow/.style={
      postaction=decorate,
      decoration={markings, mark=at position 0.5 with {\arrowreversed{angle 45}}}
    }
  ]
    % 111
    \node[p] (a10)  at (0, 0)  {};
    \node[d] (d10)  at (1, 0)  {$\cdots$};
    \node[p] (a20)  at (2, 0)  {};
    \draw[midarrow] (a10.east) -- ($(a10)!0.3!(a20)$);
    \draw[midarrow] ($(a10)!0.7!(a20)$) -- (a20.west);
    \node[p] (a30)  at (4, 0)  {};
    \node[d] (d20)  at (5,0)  {$\cdots$};
    \node[p] (a40)  at (6, 0)  {};
    \draw[midarrow] (a30.east) -- ($(a30)!0.3!(a40)$);
    \draw[midarrow] ($(a30)!0.7!(a40)$) -- (a40.west);
    \node[p] (a50)  at (8, 0)  {};
    \node[d] (d30)  at (9,0)  {$\cdots$};
    \node[p] (a60)  at (10, 0)  {};
    \draw[midarrow] (a50.east) -- ($(a50)!0.3!(a60)$);
    \draw[midarrow] ($(a50)!0.7!(a60)$) -- (a60.west);
    \node[d] (d40)  at (12,0)  {$\cdots$};

    \node[p] (a70)  at (14, 0)  {};
    \node[d] (d50)  at (15, 0)  {$\cdots$};
    \node[p] (a80)  at (16, 0)  {};
    \draw[midarrow] (a70.east) -- ($(a70)!0.3!(a80)$);
    \draw[midarrow] ($(a70)!0.7!(a80)$) -- (a80.west);
    \node[p] (a90)  at (18, 0)  {};

    % 222
    \node[p] (a11)  at (0, -1.5)  {};
    \node[d] (d11)  at (1, -1.5)  {$\cdots$};
    \node[p] (a21)  at (2, -1.5)  {};
    \draw[midarrow] (a11.east) -- ($(a11)!0.3!(a21)$);
    \draw[midarrow] ($(a11)!0.7!(a21)$) -- (a21.west);
    \node[p] (a31)  at (4, -1.5)  {};
    \node[d] (d21)  at (5, -1.5)  {$\cdots$};
    \node[p] (a41)  at (6, -1.5)  {};
    \draw[midarrow] (a31.east) -- ($(a31)!0.3!(a41)$);
    \draw[midarrow] ($(a31)!0.7!(a41)$) -- (a41.west);
    \node[p] (a51)  at (8, -1.5)  {};
    \node[d] (d31)  at (9, -1.5)  {$\cdots$};
    \node[p] (a61)  at (10, -1.5)  {};
    \draw[midarrow] (a51.east) -- ($(a51)!0.3!(a61)$);
    \draw[midarrow] ($(a51)!0.7!(a61)$) -- (a61.west);
    \node[d] (d41)  at (12,-1.5)  {$\cdots$};

    \node[p] (a71)  at (14, -1.5)  {};
    \node[d] (d51)  at (15, -1.5)  {$\cdots$};
    \node[p] (a81)  at (16, -1.5)  {};
    \draw[midarrow] (a71.east) -- ($(a71)!0.3!(a81)$);
    \draw[midarrow] ($(a71)!0.7!(a81)$) -- (a81.west);
    \node[p] (a91)  at (18, -1.5)  {};

    % ...
    \node[d] (v1)  at (1, -2.5)  {$\vdots$};
    \node[d] (v1)  at (5, -2.5)  {$\vdots$};
    \node[d] (v1)  at (9, -2.5)  {$\vdots$};
    \node[d] (v1)  at (12, -2.5)  {$\vdots$};
    \node[d] (v1)  at (15, -2.5)  {$\vdots$};
    \node[d] (v1)  at (18, -2.5)  {$\vdots$};

    % 333
    \node[p] (a12)  at (0, -3.5)  {};
    \node[d] (d12)  at (1, -3.5)  {$\cdots$};
    \node[p] (a22)  at (2, -3.5)  {};
    \draw[midarrow] (a12.east) -- ($(a12)!0.3!(a22)$);
    \draw[midarrow] ($(a12)!0.7!(a22)$) -- (a22.west);
    \node[p] (a32)  at (4, -3.5)  {};
    \node[d] (d22)  at (5, -3.5)  {$\cdots$};
    \node[p] (a42)  at (6, -3.5)  {};
    \draw[midarrow] (a32.east) -- ($(a32)!0.3!(a42)$);
    \draw[midarrow] ($(a32)!0.7!(a42)$) -- (a42.west);
    \node[p] (a52)  at (8, -3.5)  {};
    \node[d] (d32)  at (9, -3.5)  {$\cdots$};
    \node[p] (a62)  at (10, -3.5)  {};
    \draw[midarrow] (a52.east) -- ($(a52)!0.3!(a62)$);
    \draw[midarrow] ($(a52)!0.7!(a62)$) -- (a62.west);
    \node[d] (d42)  at (12,-3.5)  {$\cdots$};

    \node[p] (a72)  at (14, -3.5)  {};
    \node[d] (d52)  at (15, -3.5)  {$\cdots$};
    \node[p] (a82)  at (16, -3.5)  {};
    \draw[midarrow] (a72.east) -- ($(a72)!0.3!(a82)$);
    \draw[midarrow] ($(a72)!0.7!(a82)$) -- (a82.west);
    \node[p] (a92)  at (18, -3.5)  {};

    % lian
    \draw[midarrow] (a20.west) -- (a11.east);
    \draw[midarrow] (a21.west) -- ($(a21)!0.36!(a12)$);
    \draw[midarrow] ($(a21)!0.64!(a12)$) -- (a12.east);
    \draw[midarrow] (a22.east) -- (a30.west);

    \draw[midarrow] (a40.west) -- (a31.east);
    \draw[midarrow] (a41.west) -- ($(a41)!0.36!(a32)$);
    \draw[midarrow] ($(a41)!0.64!(a32)$) -- (a32.east);
    \draw[midarrow] (a42.east) -- (a50.west);

    \draw[midarrow] (a60.west) -- (a51.east);
    \draw[midarrow] (a61.west) -- ($(a61)!0.36!(a52)$);
    \draw[midarrow] ($(a61)!0.64!(a52)$) -- (a52.east);
    \draw[midarrow] (a62.east) -- ($(a62)!0.36!(a70)$);
    \draw[midarrow] ($(a62)!0.64!(a70)$) -- (a70.west);

    \draw[midarrow] (a80.west) -- (a71.east);
    \draw[midarrow] (a81.west) -- ($(a81)!0.36!(a72)$);
    \draw[midarrow] ($(a81)!0.64!(a72)$) -- (a72.east);

    \draw[midarrow] (a82.east) -- (a92.west);
    \draw[midarrow] (a92.north) -- ($(a92)!0.3!(a91)$);
    \draw[midarrow] ($(a92)!0.7!(a91)$)--(a91.south);
    \draw[midarrow] (a91.north) -- (a90.south);
    \draw[midarrow, rounded corners=6pt]
    (a90.north) -- (18, 0.7) -- (0,0.7) -- (a10.north);

    \draw[box] (-0.5,0.5) rectangle (2.5, -0.5);
    \draw[box] (-0.5,-1) rectangle (2.5, -2);
    \draw[box] (-0.5,-3) rectangle (2.5, -4);
    \draw[box] (3.5,0.5) rectangle (6.5, -0.5);
    \draw[box] (3.5,-1) rectangle (6.5, -2);
    \draw[box] (3.5,-3) rectangle (6.5, -4);
    \draw[box] (7.5,0.5) rectangle (10.5, -0.5);
    \draw[box] (7.5,-1) rectangle (10.5, -2);
    \draw[box] (7.5,-3) rectangle (10.5, -4);
    \draw[box] (13.5,0.5) rectangle (16.5, -0.5);
    \draw[box] (13.5,-1) rectangle (16.5, -2);
    \draw[box] (13.5,-3) rectangle (16.5, -4);
    \draw[box] (17.5,0.5) rectangle (18.5, -4);

    % sign
    \node[tlab] at (0,0.6) {$\o\d_{0,0}\setminus \lg u,v_0\rg$};
    \node[tlab] at (0,-2.6) {$\o\d_{1,0}\setminus \lg u,v_1\rg$};
    \node[tlab] at (0,-4.6) {$\o\d_{q,0}\setminus \lg u,v_q\rg$};

    \node[tlab] at (4,0.6) {$\o\d_{0,1}\setminus \lg u,v_0\rg$};
    \node[tlab] at (4,-2.6) {$\o\d_{1,1}\setminus \lg u,v_1\rg$};
    \node[tlab] at (4,-4.6) {$\o\d_{q,1}\setminus \lg u,v_q\rg$};

    \node[tlab] at (8,0.6) {$\o\d_{0,2}\setminus \lg u,v_0\rg$};
    \node[tlab] at (8,-2.6) {$\o\d_{1,2}\setminus \lg u,v_1\rg$};
    \node[tlab] at (8,-4.6) {$\o\d_{q,2}\setminus \lg u,v_q\rg$};

    \node[tlab] at (14,0.6) {$\o\d_{0,q-1}\setminus \lg u,v_0\rg$};
    \node[tlab] at (14,-2.6) {$\o\d_{1,q-1}\setminus \lg u,v_1\rg$};
    \node[tlab] at (14,-4.6) {$\o\d_{q,q-1}\setminus \lg u,v_q\rg$};
    \node[tlab] at (18,0.6) {$\o u$};

    \node[lab] at (a90.south) {$\lg u,v_0\rg$};
    \node[lab] at (a91.south) {$\lg u,v_1\rg$};
    \node[lab] at (a92) {$\lg u,v_q\rg$};
  \end{tikzpicture}
  \caption{{\footnotesize H-cycle $C$ of $\G$ in Lemma 5.3.}}
\end{figure}

\section{$M$ is a  unitary group}
By Proposition~\ref{pro:3}, we  have  three kinds of  actions for  the unitary groups,  which will be  dealt with in the following three subsections.

\subsection{Action  on isotropic points}
\begin{lem} \label{u1e}
Let $M=\SU(n,q)$ where  $n \ge 4$ is even.   acting on the set   $\O $ of  isotropic points. Then  $\G(\O, \Del_s)$  has a H-cycle.
\end{lem}
\demo Set $n=2m\ge 4$. Let $V(n,q^2)$ be the unitary  space over $\FF_{q^2}$ with the bilinear form
 $$f(X,Y)=X_1 \o Y_2^T+X_2\o Y_1^T, \,  X=(X_1,X_2), Y=(Y_1,Y_2), X_i,Y_i \in V(m,q^2).$$
  Let $\O$ be the set of isotropic points, where $|\O|=\frac{q^{2n-1}+q^n-q^{n-1}-1}{q^2-1}$.
Pick up  an  $\a \in \O$. Then
 $\D_s=\{\b \in \O \di f(\a,\b)=0\}\, \, \,  {\rm and}\, \, \,  |\D_s|=\frac{q^{2n-3}+q^n-q^{n-1}-q^2}{q^2-1}.$

Let $A$ be the matrix as given in Eq(\ref{A1}), by replacing $n$ by $m$ and $q$ by $q^2$. Set $\rho=\diag( A ,\o{A^\t})$, where $\o{A^\t}=(a_{ij}^q)_{m\times m}$ if $A^\t =(a_{ij})_{m\times m}$.
  Then $\rho \in \GU(n,q)$ and $\rho(X)=(X_1A,X_2\o{A^\t})$.  Set
\begin{equation}\label{u}\begin{array}{ll}
r=\frac{q^{2m}-1}{q^2-1}, \, &X_0=(0,1,0,\cdots,0)\in V(m,q^2),\\
X_{e,Y}=(x_1,e,x_3,\cdots,x_m)\in V(m,q^2), \, &Y=(x_1,x_3,\cdots,x_m)\in V(m-1,q^2).
\end{array}\end{equation}
Set
$$\begin{array}{llll}
B_1&=&\{\b_{1,i}\di 0\le i \le r-1\}, \quad  &\b_{1,i}=\lg(X_0A^i,0_m)\rg ;\\
B_2&=&\{\b_{2,i}\di 0\le i \le r-1\}, \quad  &\b_{2,i}=\lg(0_m,X_0\o{(A^i)^\t)}\rg ;\\
C_{Y}&=&\{\d_{Y, k}\di 0\le k \le \frac {q^{2m}-1}{q+1}-1\}, &\d_{Y, k}=\lg(X_0A^k,X_{0,Y}\o{(A^k)^\t)}\rg , \\
&&& Y=\lg \th^{q+1}\rg Y_1, Y_1\in V(m-1,q^2)\setminus \{0_{m-1}\};\\
D_{Y}&=&\{\g_{Y, k}\di 0\le k \le \frac {q^{2m}-1}{q+1}-1\}, &\g_{Y, k}=\lg(X_0A^k,X_{a, Y}\o{(A^k)^\t})\rg ,\\
&&& Y\in V(m-1,q^2),
\end{array} $$
where $a\in \FF_{q^2}^*$ is given and $a+\o a=0$.
 Check that $\O$ is a disjoint union of the above   $\lg\rho\rg$-orbits, that is
$$\O=B_1\cup B_2 \cup (\cup_{Y}  C_{Y})\cup (\cup_{Y} D_Y),$$
where the respective vectors $Y$ for $C_Y$ and $D_Y$ are given in above.
Moreover,  $\G(B_1),\G(B_2)$, $\G(\{\d_{Y,k} \di Y\})$ and $\G(\{\g_{Y,k} \di Y\})$ for any given $k$  are complete graphs.
\vskip 2mm
(1) First let $m \ge 3$.
Set $$Y_1=(1,0,0,\cdots,0),\, Y_2=(0,1,0,\cdots,0),\, Y_3=(1,1,0\cdots,0)\in V(m-1,q^2).$$

(i) Since $X_0A^{-1}\o {X_{0,Y_2}^T}= \o{X_0} \o{A} X_{0,Y_1}^T=0$, we get
$\d_{Y_1,0}\thicksim \d_{Y_2,1}=\d_{Y_2,0}^{\orho}$. In the complete graph $\G(\{\d_{Y,0} \di Y\}),$  pick a H-path:
$$U_0:\quad \d_{Y_2,0},\cdots,\d_{Y_3,0},\d_{Y_1,0},$$
and then we may get a H-path $U$ in $\G(\cup_Y C_Y)$:
$$U:\quad \d_{Y_1,0},U_0^{\orho},U_0^{\orho^2},\cdots,U_0^{\orho^{r(q-1)-1}},\d_{Y_2,0},\cdots,\d_{Y_3,0}.$$

(ii) Since  $\g_{Y_1,0}\thicksim \g_{Y_2,1}=\g_{Y_2,0}^{\orho}$.
In $\G(\{\g_{Y,0} \di Y\})$  we pick a H-path:
$$W_0:\quad \g_{Y_2,0},\g_{Y_3,0},\cdots,\g_{Y_1,0},$$
and then get  a H-path $W$ in $\G(\cup_Y D_Y)$:
$$W:\quad \g_{Y_3,0},\cdots,\g_{Y_1,0},W_0^{\orho},W_0^{\orho^2},\cdots,W_0^{\orho^{r(q-1)-1}},\g_{Y_2,0}.$$

(iii) Take  a H-path $\b_{2,1}\to\b_{2,r-1}$ in $\G(B_2)$ and a H-path  $\b_{1,r-1}\to\b_{1,1}$ in $\G(B_1)$.
Then we get a H-cycle $C$ for $\G$ (for a similar illustration, one  may refer to Figure 2):
$$C:\quad  U, \b_{2,1} \to \b_{2,r-1},W,\b_{1,r-1}\to \b_{1,1},\d_{Y_1,0}.$$

(2) Secondly, let  $m =2$.  Recall  $$X_0=(0,1),\,  A={\footnotesize\left(\begin{array}{cc} 0 &1\\a_0 &a_1\end{array}\right)}, \, A^{-1}={\footnotesize \left(\begin{array}{cc} -a_1a_0^{-1} &a_0^{-1}\\1 &0\end{array}\right)}.$$

(i)  Reset
$$\d_{Y_1,0}=\lg(0,1,-1,0)\rg , \, \d_{Y_2,0}=\lg(0,1,a_0,0)\rg , \, \d_{Y_3,0}=\lg(0,1,b,0)\rg,$$ where $\d_{Y_3,0}\ne \d_{Y_1,0},\d_{Y_2,0}$, noting $a_0 \in \FF_{q^2} \backslash \FF_q$.
Check that $\d_{Y_1,0}\thicksim \d_{Y_2,1}=\d_{Y_2,0}^{\orho}$. In $\G(\{\d_{Y,0} \di Y\})$ we chose  a H-path:
$$U_0: \quad \d_{Y_2,0},\cdots,\d_{Y_3,0},\d_{Y_1,0},$$
 and then  get a H-path $U$ in $\G(\cup_Y C_Y)$:
$$U:\quad \d_{Y_1,0},U_0^{\orho},U_0^{\orho^2},\cdots,U_0^{\orho^{r(q-1)-1}},\d_{Y_2,0},\cdots,\d_{Y_3,0}.$$

(ii) Note that $a_1\neq0$. Let
$$\g_{Y_4,0}=\lg(0,1,0,a)\rg,\, \g_{Y_5,0}=\lg(0,1,-\o aa_1,a)\rg,\, \g_{Y_6,0}=\lg(0,1,c,a)\rg , $$  where
$c\ne 0, - aa_1$.  Check that $\g_{Y_4,0}\thicksim \g_{Y_5,1}=\g_{Y_5,0}^{\orho}$. Then in  $\G(\{\g_{Y,0} \di Y\})$ we pick up  a H-path:
$$W_0:\quad \g_{Y_5,0},\cdots,\g_{Y_6,0},\g_{Y_4,0},$$
and then  get a H-path $W$ in $\G(\cup_Y D_Y)$:
$$W:\quad \g_{Y_4,0},W_0^{\orho},W_0^{\orho^2},\cdots,W_0^{\orho^{r(q-1)-1}},\g_{Y_5,0},\cdots,\g_{Y_6,0}.$$

(iii) Take a H-path $\b_{2,1}\to\b_{2,0}$ in $\G(B_2)$, and a H-path $\b_{1,r-2}\to\b_{1,0}$ in $\G(B_1)$.
Then we get a H-cycle $C$  for $\G$:
$$C:\quad U,\b_{2,1}\to \b_{2,0},(W^{-1})^{\orho^{-1}},\b_{1,r-2}\to \b_{1,0},\d_{Y_1,0}.$$
\qed

\begin{lem} \label{u1o}
Let $M=\SU(n,q)$ where  $n\ge 3$ is odd.   acting on the set   $\O $ of  isotropic points. Then  $\G(\O, \Del_s)$  has a H-cycle.
\end{lem}
\demo Set $n=2m+1\ge 3$.  Let $V=(n,q^2)$ be the $n$-dimensional unitary space over $\FF_{q^2}$ with the bilinear form
$$f(X,Y)=a\o b+X_1 \o Y_2^T+X_2\o Y_1^T, \,\, {\rm for}\,\,  X=(a,X_1,X_2),Y=(b,Y_1,Y_2), $$
where $a,b\in \FF_{q^2}, X_i,Y_i \in V(m,q^2).$
Let $\O$ be the set of isotropic points so that $|\O|=\frac{q^{2n-1}-q^n+q^{n-1}-1}{q^2-1}$. For any $\a \in \O$, we have
$\D_s=\{\b \in \O | f(\a,\b)=0\} \, \, {\rm and}\, \, |\D_s|=\frac{q^{2n-3}-q^n+q^{n-1}-q^2}{q^2-1}.$

Let $A$ be as Lemma~\ref{u1e}.(1) and   $\rho=\diag(1,A,\o{A^\t})$. Then $\rho \in \GU(n,q)$ and $\rho(X)=(a,X_1A,X_2\o{A^\t})$.  Let $r, X_0, X_{e,Y}$ have the same meaning with Eq(\ref{u}).
Set
$$\begin{array}{llll}
B_1&=&\{\b_{1,i}\di 0\le i \le r-1\}, \quad  &\b_{1,i}=\lg(0,X_0A^i,0_m)\rg ;\\
B_2&=&\{\b_{2,i}\di 0\le i \le r-1\}, \quad  &\b_{2,i}=\lg(0,0_m,X_0\o{(A^i)^\t)}\rg ;\\
C_{Y}&=&\{\d_{Y, k}\di 0\le k \le \frac {q^{2m}-1}{q+1}-1\},\quad &\d_{Y, k}=\lg(0,X_0A^k,X_{0,Y}\o{(A^k)^\t)}\rg , \\
&&& Y=\lg \th^{q+1}\rg Y_1, Y_1\in V(m-1,q^2)\setminus \{0_{m-1}\};\\
D_{Y}&=&\{\g_{Y, k}\di 0\le k \le \frac {q^{2m}-1}{q+1}-1\},\quad &\g_{Y, k}=\lg(0,X_0A^k,X_{a, Y}\o{(A^k)^\t})\rg ,\\
&&& Y\in V(m-1,q^2),a+\o a=0,\, {\rm and}\, a \ne 0;\\
E_{b,Y}&=&\{\ph_{b,Y;k} \di 0\le k \le q^{2m}-2\},\quad &\ph_{b,Y;k}=\lg(1,X_0A^k,X_{b, Y}\o{(A^k)^\t})\rg,  b\in \FF_{q^2},\\
&&&  b+\o b=-1, Y\in V(m-1,q^2).
\end{array} $$
Check that $\O$ is a disjoint union of the above   $\lg\rho\rg$-orbits, that is
$$\O=B_1\cup B_2 \cup (\cup_{Y} C_{Y})\cup (\cup_{Y} D_Y)\cup (\cup_{b,Y} E_{b,Y})$$
where the respective vectors $Y$ for $C_Y$, $D_Y$  and $E_{b,Y}$ are given in above.
Moreover,   $\G(B_1)$, $\G(B_2)$, $\G(\{\d_{Y,k} \di Y\})$,  $\G(\{\g_{Y,k} \di Y\})$  for any given $k$ and  $\G(\{\ph_{b,Y;k} \di Y\})$  for any given $k$ and $b$,  are complete graphs.
\vskip 2mm
(1) First $m \ge 3$.
Let $$Y_1=(1,0,0,\cdots,0),\, \, Y_2=(0,1,0,\cdots,0)\, \, Y_3=(1,1,0,\cdots,0)\in V(m-1,q^2).$$

(i) Check that $\d_{Y_1,0}\thicksim \d_{Y_2,1}$. In the complete graph $\G(\{\d_{Y,0}\di Y\})$, pick a H-path:
$$U_0:\quad \d_{Y_2,0},\cdots,\d_{Y_3,0},\d_{Y_1,0},$$
 and we get a H-path $U$ in $\G(\cup_Y C_Y)$:
$$U: \quad \d_{Y_1,0},U_0^{\orho},U_0^{\orho^2},\cdots,U_0^{\orho^{r(q-1)-1}},\d_{Y_2,0},\cdots,\d_{Y_3,0}.$$

(ii) Since $\g_{Y_1,0}\thicksim \g_{Y_2,1}$,  in the complete graph $\G(\g_{Y,0}\di Y)$, pick a H-path:
 $$W_0:\quad \g_{Y_2,0},\cdots,\g_{Y_3,0},\g_{Y_1,0},$$ and then get a H-path $W$ in $\G(\cup_Y D_Y)$:
$$W: \quad \g_{Y_1,0},W_0^{\orho},W_0^{\orho^2},\cdots,W_0^{\orho^{r(q-1)-1}},\g_{Y_2,0},\cdots,\g_{Y_3,0}.$$

(iii) Let $$Y_0=0_{m-1},\, \, Y_4=(b,0,0,\cdots,0),\,\, \, \,  Y_5=(0,b,0,\cdots,0)\in V(m-1,q^2).$$
 Since $1+X_0A^{-1}\o{X_{b,Y_4}^T}+\o{X_0A}X_{b,Y_5}^T=0$, we get $\ph_{b,Y_5;0}\thicksim \ph_{b,Y_4;1}=\ph_{b,Y_4;0}^\rho$.
In the complete graph $\G(\{\ph_{b,Y;0} \di Y\})$, pick a H-path:
$$R_0(b):\quad \ph_{b,Y_4;0},\ph_{b,Y_0;0}\cdots,\ph_{b,Y_5;0},$$
 and we may get a H-path $R(b)$  in $\G(\cup_{Y}E_{b,Y})$:
$$R(b): \quad \ph_{b,Y_0;0},\cdots,\ph_{b,Y_5,0},R_0(b)^\rho,R_0(b)^{\rho^2},\cdots,R_0(b)^{\rho^{q^{2m}-2}},\ph_{b,Y_4;0}.$$
Let $b_i$ where $0\le i \le q-1$ be all solutions of $\o b+b=-1$.
Since $\b_{2,1}\thicksim \ph_{b,Y_4; 0}$ and $\b_{2,1}\thicksim \ph_{b,Y_0; 2}$, we get a path (see Figure 7):
$$\begin{array}{l}R:\quad R(b_0),\b_{2,1},R(b_1)^{\rho^2},\b_{2,3},\cdots,\b_{2,2(q-1)-1},R(b_{q-1})^{\rho^{2(q-1)}}.\end{array}$$

(iv) Let $B_1'=B_1\setminus \{\b_{1,2}, \b_{1,1}\}$ and $B_2'=B_2\setminus \{\b_{2,1}, \b_{2, 3}, \cdots ,\b_{2,2(q-1)-1}\}$.
Take a H-path  $\b_{2,r-1}\to\b_{2,0}$ of $\G(B_2')$ and a  H-path  $\b_{1,2(q-1)+1}\to\b_{1,0}$ of
$\G(B_1')$.
Then we get a H-cycle $C$ for $\G$ (see Figure 8):
$$C:\quad U,\b_{2,r-1}\to\b_{2,0},(W^{-1})^{\rho},\b_{1,2},\b_{1,1},R,\b_{1,2(q-1)+1}\to\b_{1,0},\d_{Y_1,0}.$$

\begin{figure}[htbp]
 \vspace{-2mm} \centering
  \begin{tikzpicture}[
    scale=0.75, transform shape,
    p/.style={draw, circle, fill=white, inner sep=0pt, minimum size=2.6mm},
    d/.style={fill=white, inner sep=2pt},
    box/.style={draw, rounded corners=6pt},
    lab/.style={anchor=north, yshift=3pt, font=\footnotesize},
    tlab/.style={anchor=south, font=\small},
    midarrow/.style={
      postaction=decorate,
      decoration={markings, mark=at position 0.5 with {\arrow{angle 45}}}
    },
    remidarrow/.style={
      postaction=decorate,
      decoration={markings, mark=at position 0.5 with {\arrowreversed{angle 45}}}
    }
  ]
    % 111
    \node[p] (a10)  at (0, 0)  {};
    \node[p] (a20)  at (2, 0)  {};
    \node[d] (d10)  at (4, 0)  {$\cdots$};
    \node[p] (a30)  at (6, 0)  {};
    \node[p] (a40)  at (8, 0)  {};
    \node[d] (d20)  at (9, 0)  {$\cdots$};
    \node[p] (a50)  at (10, 0)  {};
    \node[p] (a60)  at (12, 0)  {};
    %\draw[midarrow] (a10.east) -- (a20.west);
    \draw[midarrow] (a20.east) -- (d10.west);
    \draw[midarrow] (d10.east) -- (a30.west);
    \draw[midarrow] (a40.east) -- ($(a40)!0.3!(a50)$);
    \draw[midarrow] ($(a40)!0.7!(a50)$) -- (a50.west);
    \draw[midarrow] (a50.east) -- (a60.west);
    % 222
    \node[p] (a11)  at (0, -1)  {};
    \node[p] (a21)  at (2, -1)  {};
    \node[d] (d11)  at (4, -1)  {$\cdots$};
    \node[p] (a31)  at (6, -1)  {};
    \node[p] (a41)  at (8, -1.2)  {};
    \node[d] (d21)  at (9, -1.2)  {$\cdots$};
    \node[p] (a51)  at (10, -1.2)  {};
    \node[p] (a61)  at (12, -1.2)  {};
    \draw[midarrow] (a11.east) -- (a21.west);
    \draw[midarrow] (a21.east) -- (d11.west);
    \draw[midarrow] (d11.east) -- (a31.west);
    \draw[midarrow] (a41.east) -- ($(a41)!0.3!(a51)$);
    \draw[midarrow] ($(a41)!0.7!(a51)$) -- (a51.west);
    \draw[midarrow] (a51.east) -- (a61.west);
    % ...
    \node[d] (v1)  at (0, -2)  {$\vdots$};
    \node[d] (v2)  at (2, -2)  {$\vdots$};
    \node[d] (v3)  at (4, -2)  {$\vdots$};
    \node[d] (v4)  at (6, -2)  {$\vdots$};
    \node[d] (v5)  at (9, -2) {$\vdots$};
    \node[d] (v5)  at (12, -2) {$\vdots$};
    % 333
    \node[p] (a12)  at (0, -3)  {};
    \node[p] (a22)  at (2, -3)  {};
    \node[d] (d12)  at (4, -3)  {$\cdots$};
    \node[p] (a32)  at (6, -3)  {};
    \node[p] (a42)  at (8, -3)  {};
    \node[d] (d22)  at (9, -3)  {$\cdots$};
    \node[p] (a52)  at (10, -3)  {};
    \draw[midarrow] (a12.east) -- (a22.west);
    \draw[midarrow] (a22.east) -- (d12.west);
    \draw[midarrow] (d12.east) -- (a32.west);
    \draw[midarrow] (a42.east) -- ($(a42)!0.3!(a52)$);
    \draw[midarrow] ($(a42)!0.7!(a52)$) -- (a52.west);

    % lian
    \draw[midarrow] (a30.west) -- (a11.east);
    \draw[midarrow] (a31.west) -- ($(a31)!0.36!(a12)$);
    \draw[midarrow] ($(a31)!0.64!(a12)$) -- (a12.east);
    \draw[midarrow, rounded corners=6pt]
    (a32.south) -- (6,-3.7) -- (-1,-3.7) -- (-1,0) -- (a10.west);
    \draw[midarrow] (a60.west) -- (a41.east);
    \draw[midarrow] (a61.west) -- ($(a61)!0.36!(a42)$);
    \draw[midarrow] ($(a61)!0.64!(a42)$) -- (a42.east);
    % block
    \draw[box] (-0.7,0.5) rectangle (0.7,-3.5);
    \draw[box] (1.3,0.5) rectangle (2.7,-3.5);
    \draw[box] (3.3,0.5) rectangle (4.7,-3.5);
    \draw[box] (5.3,0.5) rectangle (6.7,-3.5);
    \draw[box] (7.3,0.5) rectangle (10.7,-0.5);
    \draw[box] (7.3,-0.7) rectangle (10.7,-1.7);
    \draw[box] (7.3,-2.5) rectangle (10.7,-3.5);
    % sign
    \node[lab] at (a10.south) {$\ph_{b,Y_4;0}$};
    \node[lab] at (a20.south) {$\ph_{b,Y_0;0}$};
    \node[lab] at (a30.south) {$\ph_{b,Y_5;0}$};
    \node[lab] at (a11.south) {$\ph_{b,Y_4;1}$};
    \node[lab] at (a21.south) {$\ph_{b,Y_0;1}$};
    \node[lab] at (a31.south) {$\ph_{b,Y_5;1}$};
    \node[lab] at (a12.south) {$\ph_{b,Y_4;q^{2m}-2}$};
    \node[lab] at (a22.south) {$\ph_{b,Y_0;q^{2m}-2}$};
    \node[lab] at (a32.south) {$\ph_{b,Y_5;q^{2m}-2}$};
    \node[lab] at (a40.south) {$\ph_{b_0,Y_0;0}$};
    \node[lab] at (a50.south) {$\ph_{b_0,Y_4;0}$};
    \node[lab] at (a60.south) {$\b_{2,1}$};
    \node[lab] at (a41.south) {$\ph_{b_1,Y_0;2}$};
    \node[lab] at (a51.south) {$\ph_{b_1,Y_4;2}$};
    \node[lab] at (a61.south) {$\b_{2,3}$};
    \node[lab] at (a42.south) {$\ph_{b_{q-1},Y_0;2(q-1)}$};
    \node[lab] at (10.5,-3.1) {$\ph_{b_{q-1},Y_4;2(q-1)}$};

    \node[tlab] at (-1.2,0) {$R(b)$};
    \node[tlab] at (0,0.5) {$E_{b,Y_4}$};
    \node[tlab] at (2,0.5) {$E_{b,Y_0}$};
    \node[tlab] at (6,0.5) {$E_{b,Y_5}$};
    \node[tlab] at (8,0.5) {$R(b_0)$};
    \node[tlab] at (8,-2.3) {$R(b_1)$};
    \node[tlab] at (8,-4.1) {$R(b_{q-1})$};
  \end{tikzpicture}
  \caption{{\footnotesize Path $R$ in  $\G$ in Lemma 6.2.(1).}}
%\end{figure}

%\begin{figure}[htbp]
  \centering
  \begin{tikzpicture}[
    scale=0.75, transform shape,
    p/.style={draw, circle, fill=white, inner sep=0pt, minimum size=2.6mm},
    d/.style={fill=white, inner sep=2pt},
    box/.style={draw, rounded corners=6pt},
    lab/.style={anchor=north, yshift=3pt, font=\footnotesize},
    tlab/.style={anchor=south, font=\small},
    midarrow/.style={
      postaction=decorate,
      decoration={markings, mark=at position 0.5 with {\arrow{angle 45}}}
    },
    remidarrow/.style={
      postaction=decorate,
      decoration={markings, mark=at position 0.5 with {\arrowreversed{angle 45}}}
    }
  ]
  \node[p] (a10)  at (0, 0)    {};
    \node[d] (d10)  at (2, 0)  {$\cdots$};
    \node[p] (a20)  at (4, 0)  {};
    \node[p] (a30)  at (6, 0)  {};
    \node[p] (c1)   at (8, 0)  {};
    \node[d] (d)    at (8.5, 0)   {$\cdots$};
    \node[p] (c2)   at (9,0)  {};
    \node[p] (b10)  at (11, 0)  {};
    \node[p] (b20)  at (13, 0)  {};
    \node[d] (d20)  at (15, 0)  {$\cdots$};
    \node[p] (b30)  at (17, 0)  {};
    \draw[midarrow] (a10.east) -- (d10.west);
    \draw[midarrow] (d10.east) -- (a20.west);
    %\draw[midarrow] (a20.east) -- (a30.west);
    %\draw[midarrow] (b10.east) -- (b20.west);
    \draw[midarrow] (b20.east) -- (d20.west);
    \draw[midarrow] (d20.east) -- (b30.west);

    % 222
    \node[p] (a11) at (0, -1)  {};
    \node[d] (d11) at (2, -1)  {$\cdots$};
    \node[p] (a21) at (4, -1)  {};
    \node[p] (a31) at (6, -1)  {};
    \node[p] (d1)  at (8, -1.2)  {};
    \node[d] (dd)  at (8.5, -1.2) {$\cdots$};
    \node[p] (d2)  at (9, -1.2)  {};
    \node[p] (b11) at (11, -1)  {};
    \node[p] (b21) at (13, -1)  {};
    \node[d] (d21) at (15, -1)  {$\cdots$};
    \node[p] (b31) at (17, -1)  {};
    \draw[midarrow] (a11.east) -- (d11.west);
    \draw[midarrow] (d11.east) -- (a21.west);
    \draw[midarrow] (a21.east) -- (a31.west);
    \draw[midarrow] (b11.east) -- (b21.west);
    \draw[midarrow] (b21.east) -- (d21.west);
    \draw[midarrow] (d21.east) -- (b31.west);

    % ...
    \node[d] (v1) at (0, -2)   {$\vdots$};
    \node[d] (v2) at (2, -2) {$\vdots$};
    \node[d] (v3) at (4, -2) {$\vdots$};
    \node[d] (v4) at (6, -2)  {$\vdots$};
    \node[d] (v5) at (11, -2)   {$\vdots$};
    \node[d] (v6) at (13, -2) {$\vdots$};
    \node[d] (v7) at (15, -2) {$\vdots$};
    \node[d] (v8) at (17, -2)  {$\vdots$};
    \node[p] (r1) at (8,-2.7) {};
    \node[p] (r2) at (9,-2.7) {};

    % 333
    \node[p] (a12) at (0, -3)   {};
    \node[d] (d12) at (2, -3) {$\cdots$};
    \node[p] (a22) at (4, -3) {};
    \node[p] (a32) at (6, -3) {};
    \node[p] (d3) at (8, -3.8) {};
    \node[p] (d4) at (9, -3.8)   {};
    \node[p] (b12) at (11, -3) {};
    \node[p] (b22) at (13, -3) {};
    \node[d] (d22) at (15, -3) {$\cdots$};
    \node[p] (b32) at (17, -3) {};
    \draw[midarrow] (a12.east) -- (d12.west);
    \draw[midarrow] (d12.east) -- (a22.west);
    \draw[midarrow] (a22.east) -- (a32.west);
    \draw[midarrow] (b12.east) -- (b22.west);
    \draw[midarrow] (b22.east) -- (d22.west);
    \draw[midarrow] (d22.east) -- (b32.west);

    % lian
    \draw[midarrow] (a30.west)  -- (a11.east);
    \draw[remidarrow] (b10.east)  -- (b31.west);

    \draw[midarrow] (a31.west) -- ($(a31)!0.36!(a12)$);
    \draw[midarrow] ($(a31)!0.64!(a12)$) -- (a12.east);
    \draw[remidarrow] (b11.east) -- ($(b11)!0.36!(b32)$);
    \draw[remidarrow] ($(b11)!0.64!(b32)$) -- (b32.west);

    \draw[midarrow, rounded corners=6pt]
      (a32.south) -- (6, -4.1) -- (-1, -4.1) -- (-1, 0) -- (a10.west);
    \draw[remidarrow, rounded corners=6pt]
      (b12.south) -- (11, -4.1) -- (18, -4.1) -- (18, 0) -- (b30.east);

    \draw[midarrow, rounded corners=6pt]
      (a20.north) -- (4, 0.7) -- (8,0.7) -- (c1.north);
    \draw[midarrow, rounded corners=6pt]
      (c2.north) -- (9, 0.7) -- (13, 0.7) -- (b20.north);
    \draw[midarrow] (b10.west) -- (d4.east);
    \draw[midarrow] (d4.west) -- (d3.east);
    \draw[midarrow] (d3.east) -- (r2.west);
    \draw[midarrow] (r1.east) -- (d2.west);
    \draw[midarrow] (d1.west) -- (a30.east);

    % block
    \draw[box] (-0.7, 0.5) rectangle (0.7, -3.9);
    \draw[box] (1.3, 0.5) rectangle (2.7, -3.9);
    \draw[box] (3.3, 0.5) rectangle (4.7, -3.9);
    \draw[box] (5.3, 0.5) rectangle (6.7, -3.9);
    \draw[box] (7.3, 0.5) rectangle (9.7, -0.5);
    \draw[box] (7.3, -0.7) rectangle (9.7, -1.7);
    \draw[box] (7.3, -2.2) rectangle (9.7, -3.2);
    \draw[box] (10.3, 0.5) rectangle (11.7, -3.9);
    \draw[box] (12.3, 0.5) rectangle (13.7, -3.9);
    \draw[box] (14.3, 0.5) rectangle (15.7, -3.9);
    \draw[box] (16.3, 0.5) rectangle (17.7, -3.9);

    % sign
    \node[lab] at (a10.south)  {$\d_{Y_2,0}$};
    \node[lab] at (a20.south)  {$\d_{Y_3,0}$};
    \node[lab] at (a30.south)  {$\d_{Y_1,0}$};
    \node[lab] at (b10.south)  {$\g_{Y_1,1}$};
    \node[lab] at (b20.south)  {$\g_{Y_3,1}$};
    \node[lab] at (b30.south)  {$\g_{Y_2,1}$};

    \node[lab] at (a11.south)  {$\d_{Y_2,1}$};
    \node[lab] at (a21.south)  {$\d_{Y_3,1}$};
    \node[lab] at (a31.south)  {$\d_{Y_1,1}$};
    \node[lab] at (b11.south)  {$\g_{Y_1,2}$};
    \node[lab] at (b21.south)  {$\g_{Y_3,2}$};
    \node[lab] at (b31.south)  {$\g_{Y_2,2}$};

    \node[lab] at (a12.south)  {$\d_{Y_2,r(q-1)-1}$};
    \node[lab] at (a22.south)  {$\d_{Y_3,r(q-1)-1}$};
    \node[lab] at (a32.south)  {$\d_{Y_1,r(q-1)-1}$};
    \node[lab] at (b12.south)  {$\g_{Y_1,0}$};
    \node[lab] at (b22.south)  {$\g_{Y_3,0}$};
    \node[lab] at (b32.south)  {$\g_{Y_2,0}$};

    \node[lab] at (c1.south) {$\b_{2,r-1}$};
    \node[lab] at (c2.south) {$\b_{2,0}$};
    \node[lab] at (d1.south) {$\b_{1,0}$};
    \node[lab] at (9,-0.6) {$\b_{1,2(q-1)+1}$};
    \node[lab] at (d3.south) {$\b_{1,1}$};
    \node[lab] at (d4.south) {$\b_{1,2}$};
    \node[lab] at (r1.south) {$\ph_{b_{q-1},Y_4;2(q-1)}$};
    \node[lab] at (9,-2.2) {$\ph_{b_0,Y_0;0}$};

    \node[tlab] at (0, 0.6) {$C_{Y_2}$};
    \node[tlab] at (2, 0.6) {$\cdots$};
    \node[tlab] at (4, 0.6) {$C_{Y_3}$};
    \node[tlab] at (6, 0.6) {$C_{Y_1}$};
    \node[tlab] at (11, 0.6) {$D_{Y_1}$};
    \node[tlab] at (13, 0.6) {$D_{Y_3}$};
    \node[tlab] at (15, 0.6) {$\cdots$};
    \node[tlab] at (17, 0.6) {$C_{Y_2}$};
    \node[tlab] at (8, 0.6) {$B_1'$};
    \node[tlab] at (8, -2.3) {$B_2'$};
    \node[tlab] at (8, -3.7) {$R$};

  \end{tikzpicture}

  \caption{{\footnotesize H-cycle $C$ for $\G$ in Lemma 6.2.(1).}}
\end{figure}

(2) Secondly, let $m=2$. Recall $a_0$ and  $a_1$  are two entries of the matrix $A$.
\vskip 2mm
(i) Let $$\d_{Y_1,0}=\lg(0,0,1,-1,0)\rg , \, \d_{Y_2,0}=\lg(0,0,1,a_0,0)\rg , \, \,  \d_{Y_3,0}=\lg(0,0,1,c,0)\rg ,$$
where  $\d_{Y_3,0}\ne \d_{Y_1,0}, \d_{Y_2,0}$. Then we  have  a H-path $U_0$  in $\G(\{\d_{Y,0}\di Y\})$:
$$U_0:\quad \d_{Y_2,0},\cdots,\d_{Y_3,0},\d_{Y_1,0},$$  and a H-path in $\G(\cup_Y C_Y)$:
$$U: \quad \d_{Y_1,0},U_0^{\orho},U_0^{\orho^2},\cdots,U_0^{\orho^{r(q-1)-1}},\d_{Y_2,0},\cdots,\d_{Y_3,0}.$$

(ii) Let
$$\g_{Y_4,0}=\lg(0,0,1,0, a)\rg,\, \, \g_{Y_5,0}=\lg(0,0,1,-\o aa_1,a),\rg \, \, \, \g_{Y_6,0}=\lg(0,0,1,d,a)\rg ,$$
where $d\ne 0,  -\o aa_1$.
 Then  we get a  H-path $W_0$  in $\G(\{\g_{Y,0} \di Y\})$:
 $$W_0:\quad g_{Y_5,0},\cdots,\g_{Y_6,0},\g_{Y_4,0},$$ and  a H-path $W$ in $\G(\cup_Y D_Y)$:
$$W:\quad \g_{Y_4,0},W_0^{\orho},W_0^{\orho^2},\cdots,W_0^{\orho^{r(q-1)-1}},\g_{Y_5,0},\cdots,\g_{Y_6,0}.$$

(iii) No loss, among of $q$ solutions $b_i$ of $\o b+b=-1$,  assume $-1-a_1\o b_j\neq 0$ when $1\le j \le q-1,$ while $-1-a_1\o b_0$ might be 0.
    Set
$$\begin{array}{ll} \ph_{b,Y_7;0}=\lg(1,0,1,0,b)\rg, \quad  &\ph_{b,Y_8;0}=\lg(1,0,1,-1-a_1\o b,b)\rg ,\\
 \ph_{b_0,Y_9;0}=\lg(1,0,1,-1-a_1\o b_1,b_0)\rg, &\\
 \ph_{b_1,Y_{10};0}=\lg(1,0,1,\o a_0^{-1}(-1-a_1 \o b_1),b_1)\rg, \quad &\ph_{b_0,Y_{11};0}=\lg(1,0,1,b_1\o a_1-a_1\o b_1,b_0)\rg .\end{array}$$
Then check that
$$\ph_{b_1,Y_7;0}\thicksim \ph_{b_0,Y_9;1}=\ph_{b_0,Y_9;0}^\rho,   \ph_{b_1,Y_{10};0}\thicksim \ph_{b_0,Y_{11};1},
   \ph_{b_0,Y_9;0}\neq \ph_{b_0,Y_{11};0},   \ph_{b_1,Y_{10};0}\neq \ph_{b_1,Y_8;0}.$$

   First,  we get a H-cycle in $\G(\{\ph_{b_0,Y;1}\di Y\}\cup\{\ph_{b_1,Y;0}\di Y\})$:
$$\ph_{b_1,Y_7;0},\ph_{b_0,Y_9;1}\to\ph_{b_0,Y_{11};1},\ph_{b_1,Y_{10};0}\to\ph_{b_1,Y_8;0},\ph_{b_1,Y_7;0},$$
and then   pick a H-path $R_0(b_0, b_1)$ leading $\ph_{b_1,Y_8;0}$ to $\ph_{b_1,Y_7;0}$ in this cycle.  Therefore, we get a H-cycle in $\G((\cup_Y E_{b_0,Y})\cup  (\cup_Y E_{b_1,Y}))$ (see Figure 9):
$$R(b_0, b_1):\quad R_0(b_0, b_1),R_0((b_0, b_1)^\rho,R_0(b_0, b_1)^{\rho^2},\cdots,R_0(b_0, b_1)^{\rho^{q^{2m}-2}},\ph_{b_1,Y_8;0},$$
so that we may  pick a H-path $\ph_{b_1,Y_7;0}\to\ph_{b_1,Y_8;1}=\ph_{b_1,Y_8;0}^\rho$ in $R(b_0,b_1)$.

 Second,   for  $2\le j\le q-1$,  in the complete graphs $\G(\{\ph_{b_j,Y;0} \di Y\})$, we may pick a H-path $R_0(b_j)$ leading $\ph_{b_j,Y_8;0}$ to $\ph_{b_j,Y_7;0}$.
Then get a H-cycle in $\G(\cup_Y E_{b_j,Y})$:
$$R(b_j):\quad R_0(b_j),R_0(b_j)^\rho,R_0(b_j)^{\rho^2},\cdots,R_0(b_j)^{\rho^{q^{2m}-2}},\ph_{b_j,Y_8;0},$$
so that we may  pick a H-path $\ph_{b_j,Y_7;0}\to\ph_{b_j,Y_8;1}=\ph_{b_j,Y_8;0}^\rho$ in $R(b_j)$.

 Third, observing  $\b_{1,1}\thicksim \ph_{b_j,Y_7;2},\b_{1,1}\thicksim \b_{2,2},\b_{2,2}\thicksim \ph_{b_j,Y_8;1}$,  we get a path (see Figure 10):
$$\begin{matrix}R:\quad \ph_{b_1,Y_7;0}\to\ph_{b_1,Y_8;1},\b_{2,2},\b_{1,1},\ph_{b_2,Y_7;2}\to\ph_{b_2,Y_8;3},\cdots,\b_{2,2(q-2)},\\
\b_{1,2(q-2)-1},\ph_{b_{q-1},Y_7;2(q-2)}\to\ph_{b_{q-1},Y_8;2(q-2)+1}.\end{matrix}$$
\begin{figure}[htbp]
  \centering
  \begin{tikzpicture}[
    scale=0.75, transform shape,
    p/.style={draw, circle, fill=white, inner sep=0pt, minimum size=2.6mm},
    d/.style={fill=white, inner sep=2pt},
    box/.style={draw, rounded corners=6pt},
    lab/.style={anchor=north, yshift=3pt, font=\footnotesize},
    tlab/.style={anchor=south, font=\small},
    midarrow/.style={
      postaction=decorate,
      decoration={markings, mark=at position 0.5 with {\arrow{angle 45}}}
    },
    remidarrow/.style={
      postaction=decorate,
      decoration={markings, mark=at position 0.5 with {\arrowreversed{angle 45}}}
    }
  ]
    % 111
    \node[p] (a10)  at (0, 0)  {};
    \node[p] (a20)  at (2, 0)  {};
    \node[d] (d10)  at (4, 0)  {$\cdots$};
    \node[p] (a30)  at (6, 0)  {};
    \node[p] (a40)  at (8, 0)  {};
    \node[d] (d20)  at (10, 0)  {$\cdots$};
    \node[p] (a50)  at (12, 0)  {};
    \draw[midarrow] (a10.east) -- (a20.west);
    \draw[midarrow] (a20.east) -- (d10.west);
    \draw[midarrow] (d10.east) -- (a30.west);
    \draw[midarrow] (a30.east) -- (a40.west);
    \draw[midarrow] (a40.east) -- (d20.west);
    \draw[midarrow] (d20.east) -- (a50.west);
    % 222
    \node[p] (a11)  at (0, -1)  {};
    \node[p] (a21)  at (2, -1)  {};
    \node[d] (d11)  at (4, -1)  {$\cdots$};
    \node[p] (a31)  at (6, -1)  {};
    \node[p] (a41)  at (8, -1)  {};
    \node[d] (d21)  at (10, -1)  {$\cdots$};
    \node[p] (a51)  at (12, -1)  {};
    \draw[midarrow] (a11.east) -- (a21.west);
    \draw[midarrow] (a21.east) -- (d11.west);
    \draw[midarrow] (d11.east) -- (a31.west);
    \draw[midarrow] (a31.east) -- (a41.west);
    \draw[midarrow] (a41.east) -- (d21.west);
    \draw[midarrow] (d21.east) -- (a51.west);
    % ...
    \node[d] (v1)  at (0, -2)  {$\vdots$};
    \node[d] (v2)  at (2, -2)  {$\vdots$};
    \node[d] (v3)  at (4, -2)  {$\vdots$};
    \node[d] (v4)  at (6, -2)  {$\vdots$};
    \node[d] (v5)  at (8, -2) {$\vdots$};
    \node[d] (v6)  at (10, -2) {$\vdots$};
    \node[d] (v7)  at (12, -2) {$\vdots$};
    % 333
    \node[p] (a12)  at (0, -3)  {};
    \node[p] (a22)  at (2, -3)  {};
    \node[d] (d12)  at (4, -3)  {$\cdots$};
    \node[p] (a32)  at (6, -3)  {};
    \node[p] (a42)  at (8, -3)  {};
    \node[d] (d22)  at (10, -3)  {$\cdots$};
    \node[p] (a52)  at (12, -3)  {};
    \draw[midarrow] (a12.east) -- (a22.west);
    \draw[midarrow] (a22.east) -- (d12.west);
    \draw[midarrow] (d12.east) -- (a32.west);
    \draw[midarrow] (a32.east) -- (a42.west);
    \draw[midarrow] (a42.east) -- (d22.west);
    \draw[midarrow] (d22.east) -- (a52.west);

    % lian
    %\draw[midarrow] (a10.east) -- (a51.west);
    \draw[remidarrow] (a11.east) -- ($(a11)!0.36!(a52)$);
    \draw[remidarrow] ($(a11)!0.64!(a52)$) -- (a52.west);
    \draw[remidarrow, rounded corners=6pt]
    (a12.south) -- (0,-3.7) -- (12.9,-3.7) -- (12.9,0) -- (a50.east);
    % block
    \draw[box] (-0.7,0.5) rectangle (0.7,-3.5);
    \draw[box] (1.3,0.5) rectangle (2.7,-3.5);
    \draw[box] (3.3,0.5) rectangle (4.7,-3.5);
    \draw[box] (5.3,0.5) rectangle (6.7,-3.5);
    \draw[box] (7.3,0.5) rectangle (8.7,-3.5);
    \draw[box] (9.3,0.5) rectangle (10.7,-3.5);
    \draw[box] (11.3,0.5) rectangle (12.7,-3.5);

    % sign
    \node[lab] at (a10.south) {$\ph_{b_1,Y_7;0}$};
    \node[lab] at (a20.south) {$\ph_{b_0,Y_9;1}$};
    \node[lab] at (a30.south) {$\ph_{b_0,Y_{11};1}$};
    \node[lab] at (a40.south) {$\ph_{b_1,Y_{10};0}$};
    \node[lab] at (a50.south) {$\ph_{b_1,Y_8;0}$};
    \node[lab] at (a11.south) {$\ph_{b_1,Y_7;1}$};
    \node[lab] at (a21.south) {$\ph_{b_0,Y_9;2}$};
    \node[lab] at (a31.south) {$\ph_{b_0,Y_{11};2}$};
    \node[lab] at (a41.south) {$\ph_{b_1,Y_{10};1}$};
    \node[lab] at (a51.south) {$\ph_{b_1,Y_8;1}$};
    \node[lab] at (a12.south) {$\ph_{b_1,Y_7;q^{2m}-2}$};
    \node[lab] at (a22.south) {$\ph_{b_0,Y_9;0}$};
    \node[lab] at (a32.south) {$\ph_{b_0,Y_{11};0}$};
    \node[lab] at (a42.south) {$\ph_{b_1,Y_{10};q^{2m}-2}$};
    \node[lab] at (a52.south) {$\ph_{b_1,Y_8;q^{2m}-2}$};

    \node[tlab] at (-1.4,0) {$R(b_0,b_1)$};
    \node[tlab] at (0,0.5) {$E_{b_1,Y_7}$};
    \node[tlab] at (2,0.5) {$E_{b_0,Y_9}$};
    \node[tlab] at (6,0.5) {$E_{b_0,Y_{11}}$};
    \node[tlab] at (8,0.5) {$E_{b_1,Y_{10}}$};
    \node[tlab] at (12,0.5) {$E_{b_1,Y_8}$};

  \end{tikzpicture}
  \caption{{\footnotesize Path $R(b_0,b_1)$ for $\G$ in Lemma 6.2.(2).}}
%\end{figure}
\vskip 2mm
%\begin{figure}[htbp]
  \centering
  \begin{tikzpicture}[
    scale=0.75, transform shape,
    p/.style={draw, circle, fill=white, inner sep=0pt, minimum size=2.6mm},
    d/.style={fill=white, inner sep=2pt},
    box/.style={draw, rounded corners=6pt},
    lab/.style={anchor=north, yshift=3pt, font=\footnotesize},
    tlab/.style={anchor=south, font=\small},
    midarrow/.style={
      postaction=decorate,
      decoration={markings, mark=at position 0.5 with {\arrow{angle 45}}}
    },
    remidarrow/.style={
      postaction=decorate,
      decoration={markings, mark=at position 0.5 with {\arrowreversed{angle 45}}}
    }
  ]
    % 111
    \node[p] (a10)  at (0, 0)  {};
    \node[d] (d10)  at (2, 0)  {$\cdots$};
    \node[p] (a20)  at (4, 0)  {};
    \node[p] (a30)  at (6, 0)  {};
    \node[d] (d20)  at (7, 0)  {$\cdots$};
    \node[p] (a40)  at (8, 0)  {};
    \node[p] (a50)  at (10, 0)  {};
    \node[p] (a60)  at (12, 0)  {};
    \draw[midarrow] (a10.east) -- (d10.west);
    \draw[midarrow] (d10.east) -- (a20.west);
    \draw[midarrow] (a30.east) -- ($(a30)!0.3!(a40)$);
    \draw[midarrow] ($(a30)!0.7!(a40)$) -- (a40.west);
    \draw[midarrow] (a40.east) -- (a50.west);
    \draw[midarrow] (a50.east) -- (a60.west);
    % 222
    \node[p] (a11)  at (0, -1)  {};
    \node[d] (d11)  at (2, -1)  {$\cdots$};
    \node[p] (a21)  at (4, -1)  {};
    \node[p] (a31)  at (6, -1.2)  {};
    \node[d] (d21)  at (7, -1.2)  {$\cdots$};
    \node[p] (a41)  at (8, -1.2)  {};
    \node[p] (a51)  at (10, -1.2)  {};
    \node[p] (a61)  at (12, -1.2)  {};
    \draw[midarrow] (a11.east) -- (d11.west);
    \draw[midarrow] (d11.east) -- (a21.west);
    \draw[midarrow] (a31.east) -- ($(a31)!0.3!(a41)$);
    \draw[midarrow] ($(a31)!0.7!(a41)$) -- (a41.west);
    \draw[midarrow] (a41.east) -- (a51.west);
    \draw[midarrow] (a51.east) -- (a61.west);
    % ...
    \node[d] (v1)  at (0, -2)  {$\vdots$};
    \node[d] (v2)  at (2, -2)  {$\vdots$};
    \node[d] (v3)  at (4, -2)  {$\vdots$};
    \node[d] (v4)  at (7, -2)  {$\vdots$};
    \node[d] (v5)  at (10, -2)  {$\vdots$};
    \node[d] (v6)  at (12, -2)  {$\vdots$};

    % 333
    \node[p] (a12)  at (0, -3)  {};
    \node[d] (d12)  at (2, -3)  {$\cdots$};
    \node[p] (a22)  at (4, -3)  {};
    \node[p] (a32)  at (6, -3)  {};
    \node[d] (d22)  at (7, -3)  {$\cdots$};
    \node[p] (a42)  at (8, -3)  {};
    \draw[midarrow] (a12.east) -- (d12.west);
    \draw[midarrow] (d12.east) -- (a22.west);
    \draw[midarrow] (a32.east) -- ($(a32)!0.3!(a42)$);
    \draw[midarrow] ($(a32)!0.7!(a42)$) -- (a42.west);

    % lian
    %\draw[remidarrow] (a10.east) -- (a21.west);
    \draw[remidarrow] (a11.west) -- ($(a11)!0.36!(a22)$);
    \draw[remidarrow] ($(a11)!0.64!(a22)$) -- (a22.east);
    \draw[remidarrow, rounded corners=6pt]
    (a12.south) -- (0,-3.7) -- (5,-3.7) -- (5,0) -- (a20.east);
    \draw[midarrow] (a60.west) -- (a31.east);
    \draw[midarrow] (a61.west) -- ($(a61)!0.36!(a32)$);
    \draw[midarrow] ($(a61)!0.64!(a32)$) -- (a32.east);
    % block
    \draw[box] (-0.7,0.5) rectangle (0.7,-3.5);
    \draw[box] (1.3,0.5) rectangle (2.7,-3.5);
    \draw[box] (3.3,0.5) rectangle (4.7,-3.5);
    \draw[box] (5.3,0.5) rectangle (8.7,-0.5);
    \draw[box] (5.3,-0.7) rectangle (8.7,-1.7);
    \draw[box] (5.3,-2.5) rectangle (8.7,-3.5);
    % sign
    \node[lab] at (a10.south) {$\ph_{b,Y_7;0}$};
    \node[lab] at (a20.south) {$\ph_{b,Y_8;0}$};
    \node[lab] at (a11.south) {$\ph_{b,Y_7;1}$};
    \node[lab] at (a21.south) {$\ph_{b,Y_8;1}$};
    \node[lab] at (a12.south) {$\ph_{b,Y_7;q^{2m}-2}$};
    \node[lab] at (a22.south) {$\ph_{b,Y_8;q^{2m}-2}$};
    \node[lab] at (a30.south) {$\ph_{b_1,Y_7;0}$};
    \node[lab] at (a40.south) {$\ph_{b_1,Y_8;1}$};
    \node[lab] at (a50.south) {$\b_{2,2}$};
    \node[lab] at (a60.south) {$\b_{1,1}$};
    \node[lab] at (a31.south) {$\ph_{b_1,Y_7;2}$};
    \node[lab] at (a41.south) {$\ph_{b_1,Y_8;3}$};
    \node[lab] at (a51.south) {$\b_{2,4}$};
    \node[lab] at (a61.south) {$\b_{1,3}$};
    \node[lab] at (a32.south) {$\ph_{b_{q-1},Y_7;2(q-2)}$};
    \node[lab] at (8.4,-3.2) {$\ph_{b_{q-1},Y_8;2(q-2)+1}$};

    \node[tlab] at (-1.2,0) {$R(b)$};
    \node[tlab] at (0,0.5) {$E_{b,Y_7}$};
    \node[tlab] at (2,0.5) {$\cdots$};
    \node[tlab] at (4,0.5) {$E_{b,Y_8}$};
    \node[tlab] at (6,0.5) {$R(b_1,b_0)$};
    \node[tlab] at (6,-2.3) {$R(b_2)$};
    \node[tlab] at (6,-4.1) {$R(b_{q-1})$};

  \end{tikzpicture}
  \caption{{\footnotesize Path $R$ for $\G$ in Lemma 6.2.(2).}}
\end{figure}

(iv) Take a H-path  $\b_{2, 1}\to\b_{2,0}$ of $\G(B_2\setminus \{\b_{2,2}, \b_{2, 4}, \cdots ,\b_{2,2(q-1)}\})$
and a  H-path  $\b_{1,2(q-1)-1}\to\b_{1,0}$  of
$\G(B_1\setminus (\{\b_{1,1}, \b_{1,3}, \cdots ,\b_{1, 2(q-2)-1} \}\cup \{\b_{1, r-2}, \b_{1,r-1}\} )$.
 Finally,  we get a H-cycle $C$ for $\G$:
$$C:\quad U,\b_{2,1}\to\b_{2,0},(W^{-1})^{\rho^{-1}},\b_{1,r-2},\b_{1,r-1},R,
\b_{2,2(q-1)},\b_{1,2(q-1)-1}\to\b_{1,0},\d_{Y_1,0}.$$
\qed

\subsection{Action  on nonisotropic points}
\begin{lem}\label{u3e}
Let $M=\SU(n,2)$ where $n$ is even, acting on the set $\O$ of   nonisotropic points. Then $\G(\O, \Del_s)$  has a H-cycle.
\end{lem}
\demo Set $n=2m$. Let $f$ and $\rho$  be as in Lemma~\ref{u1e}.  We  adapt the notations of Eq(~\ref{u}). Let
 $\O$  the set of nonisotropic points so that $|\O|=\frac{2^{2n-1}-2^{n-1}}{3}$. Note that $\FF_{2^2}^*=\lg\th\rg$. Pick up an $\a\in\O$, then we have
 $\D_s=\{\b\in\O \di f(\a,\b)=0\}\, \, \, {\rm and}\, \, \, |\D_s|=\frac{2^{2n-3}+2^{n-2}}{3}.$
 Set
$$\begin{array}{llll}
C_{Y}&=&\{\d_{Y, k}\di 0\le k \le \frac {2^{2m}-1}{3}-1\},\quad &\d_{Y, k}=\lg(X_0A^k,X_{\th,Y}\o{(A^k)^\t)}\rg , \\
&&& Y\in V(m-1,2^2);\\
D_{Y}&=&\{\g_{Y, k}\di 0\le k \le \frac {2^{2m}-1}{3}-1\},\quad &\g_{Y, k}=\lg(X_0A^k,X_{\o\th, Y}\o{(A^k)^\t})\rg ,\\
&&& Y\in V(m-1,2^2).\\
\end{array} $$
Check that $\O$ is a disjoint union of the above   $\lg\rho\rg$-orbits, that is
$$\O=(\cup_{Y\in V(m-1, 4)} C_{Y})\cup (\cup_{Y\in V(m-1, 4)} D_Y).$$
Moreover, $\G(\{\d_{Y,k}\}\cup\{\g_{Y',k}\} \di Y,Y')$ is a complete bipartite graph for any  given $k$.
Let
$$\begin{array}{lll} Y_1=(1,0,0,\cdots,0), \quad &Y_2=(0,1,0,\cdots,0)\in V(m-1,2^2)\quad {\rm for}\, m\ge 3;\\
 \d_{Y_1,0}=\lg(0,1,0,\th)\rg , \quad &\g_{Y_2,0}=\lg(0,1,-\o\th a_1,\o\th)\rg  \quad  {\rm for}\, m=2.\end{array}$$
Then we get $\d_{Y_1,0}\thicksim\g_{Y_2,1}=\g_{Y_2,0}^{\orho}$.
In the complete bipartite graph $\G(\{\d_{Y,0}\}\cup\{\g_{Y',0}\} \di Y,Y')$, pick a H-path $U$ leading $\g_{Y_2,0}$ to $\d_{Y_1,0}$. Then we get a H-cycle $C$  for $\G$:
$$ C:\quad U,U^{\orho},U^{\orho^2},\cdots,U^{\orho^{\frac{2^{2m}-1}{3}-1}},\g_{Y_2,0}.$$
\qed

\begin{lem}\label{u3o}
Let $M=\SU(n,2)$ where $n\ge 5$ is odd, acting on the set $\O$ of   nonisotropic points. Then $\G(\O, \Del_s)$  has a H-cycle.
\end{lem}
\demo Set   $n=2m+1\ge 5$. Let $f$  and $\rho$  be as in Lemma~\ref{u1o}.  Adapt the notations of Eq(\ref{u}).  Let $\O$ be the set of nonisotropic points so that $|\O|=\frac{2^{2n-1}+2^{n-1}}{3}$.
For any $\a\in\O$, we have
$\D_s=\{\b\in\O \di f(\a,\b)=0\}\, \, \, {\rm and}\, \, \, |\D_s|=\frac{2^{2n-3}-2^{n-2}}{3}.$ Let $\G=\G(\O, \Del_s)$. Then
 $\{\a_1, \a_2\}\in E(\G)$ if and only if  $f(\a_1, \a_2)=0$.

Set
$$\begin{array}{llll}
\a &=&\lg(1,0_m,0_m)\rg;\\
B_1&=&\{\b_{1,k} \di 0\le k \le 2^{2m}-2\},\quad &\b_{1,k}=\lg(1,X_0A^k,0_m)\rg;\\
B_2&=&\{\b_{2,k} \di 0\le k \le 2^{2m}-2\},\quad &\b_{2,k}=\lg(1,0_m,X_0\o{(A^k)^\t})\rg;\\
C_{Y}&=&\{\d_{Y, k}\di 0\le k \le \frac {2^{2m}-1}{3}-1\},\quad &\d_{Y, k}=\lg(0,X_0A^k,X_{\th,Y}\o{(A^k)^\t)}\rg ;\\
&&& Y\in V(m-1,2^2);\\
D_{Y}&=&\{\g_{Y, k}\di 0\le k \le \frac {2^{2m}-1}{3}-1\},\quad &\g_{Y, k}=\lg(0,X_0A^k,X_{\o\th, Y}\o{(A^k)^\t})\rg ,\\
&&& Y\in V(m-1,2^2);\\
E_{0,Y}&=&\{\ph_{0,Y;k} \di 0\le k \le 2^{2m}-2\},\quad &\ph_{0,Y;k}=\lg(1,X_0A^k,X_{0,Y}\o{(A^k)^\t})\rg, \\
&&& Y\in V(m-1,2^2)\setminus \{0_{m-1}\};\\
E_{1,Y}&=&\{\ph_{1,Y,k} \di 0\le k \le 2^{2m}-2\},\quad &\ph_{1,Y;k}=\lg(1,X_0A^k,X_{1,Y}\o{(A^k)^\t})\rg ,\\
&&& Y\in V(m-1,2^2).
\end{array} $$
Check that $\O$ is a disjoint union of the above   $\lg\rho\rg$-orbits, that is
$$\O=\{ \a\}\cup B_1\cup B_2\cup (\cup_{Y} C_{Y})\cup (\cup_Y D_Y) \cup (\cup _Y E_{0,Y}) \cup (\cup _Y E_{1,Y}),$$
where the respective vectors $Y$ for $C_Y$, $D_Y$, $E_{0,Y}$ and  $E_{1,Y}$ are given in above.
Then $\G(\{\b_{1,k},\ph_{0,Y;k} \di Y\}\cup\{\ph_{1,Y';k} \di Y'\})$ is a complete bipartite graph for any given $k$, and every vertex  in $C_Y$ and $D_Y$ coincides with one of that   in  $C_Y$ and $D_Y$ in Lemma~\ref{u3e}, by adding 0 in the first coordinate.
  Note that  $\a$ is adjacent to all vertices of $C_Y$ and $D_Y$.

\vskip 2mm
(1) First, let $m\ge3$.
Let $$Y_0=(0,0,0,\cdots,0)\quad {\rm and}\quad  Y_1=(1,0,0,\cdots,0)\in V(m-1,2^2).$$
\vskip 2mm
(i) With the way as Lemma~\ref{u3e},  take a H-path $\d_{Y_0,1}\to\g_{Y_0,1}$ in $\G((\cup_Y C_Y)\cup (\cup_Y D_Y))$ (see Figure 11).
Since $\a$ is adjacent to all vertices of  $\G((\cup_Y C_Y)\cup (\cup_Y D_Y))$,
 Then  we get a  H-path   $U$  of  $\G((\cup_Y C_Y)\cup (\cup_Y D_Y)\cup \{\a\})$:
$$U:\quad \d_{Y_0,1}, \a\to\g_{Y_0,1}$$

\vskip 2mm
(ii)  Check that $\ph_{1,Y_0;0}\thicksim\ph_{1,Y_1;1}=\ph_{1,Y_1;0}^\rho$.
 Since $\b_{2,k}$ is adjacent to all vertices of the complete bipartite graph $\G(\{\b_{1,k},\ph_{0,Y;k}\di Y\}\cup\{\ph_{1,Y';k}\di Y'\})$ for given $k$,
one may take a  H-path $W_0$:  $\ph_{1,Y_1;0}\to \ph_{1,Y_0;0}$ of
 $\G(\{\b_{1,0},\ph_{0,Y;0}\di Y\}\cup\{\b_{2,0}, \ph_{1,Y';0}\di Y'\})$
 and then   get a H-cycle in $\G((\cup_YE_{0,Y})\cup (\cup_{Y} E_{1,Y})$  $\cup B_1\cup B_2)$ (see Figure 12):
$$W:\quad W_0,W_0^\rho,W_0^{\rho^2},\cdots,W_0^{\rho^{2^{2m}-2}},\ph_{1,Y_1;0}.$$

\vskip 2mm
(iii)  Check  $\b_{1,0}\thicksim \d_{Y_0,1}$ and $\b_{2,0}\thicksim \g_{Y_0,1}$. Pick a H-path $\b_{2,0}\to\b_{1,0}$ in $W$, and
  Finally,   we get a H-cycle $C$ for $\G$:
$$C:\quad \b_{2,0}\to\b_{1,0},\d_{Y_0,1}, \a \to\g_{Y_0,1},\b_{2,0}.$$

(2) Secondly, let $m=2.$
Now $X_0=(0,1)$. Let $\ph_{0,Y_3;0}=\lg(1,0,1,1,0)\rg$. Check that $\b_{1,0}\thicksim \ph_{0,Y_3,1}$.  Pick a H-path $W_0:$ $\ph_{0,Y_3;0}\to\b_{1,0}$ of   $\G(\{\b_{1,0},\b_{2,0},$ $\ph_{0,Y;0}\di Y\}\cup\{\ph_{1,Y';0}\di Y'\})$, and then
  get a H-cycle $W$  of  $\G((\cup_Y E_{0,Y})\cup(\cup_Y E_{1,Y})\cup B_1\cup B_2)$:
$$W:\quad W_0,W_0^\rho,W_0^{\rho^2},\cdots,W_0^{\rho^{2^{2m}-2}},\ph_{0,Y_3;0}.$$
Let $$\d_{Y_1,0}=\lg(0,0,1,0,\th)\rg \quad {\rm and}\quad  \g_{Y_1,0}=\lg(0,0,1,0,\o\th)\rg.$$
With the way as Lemma~\ref{u3e}, pick a H-path $\d_{Y_1,1}\to\g_{Y_1,1}$ in $\G((\cup_YC_Y)\cup(\cup_Y D_Y))$. Let $\ph_{1,Y_4;0}=\lg(1,0,1,\o a_0^{-1}\o a_1,1)\rg$. Check  $\b_{1,0}\thicksim\d_{Y_1,1}$ and $\ph_{1,Y_4,0}\thicksim \g_{Y_1,1}$.
Then we have a H-path $\d_{Y_1,1}, \a\to\g_{Y_1,1}$ of  $\G((\cup_YC_Y)\cup(\cup_Y D_Y)\cup \{\a\})$.
 Pick a H-path $\ph_{1,Y_4;0}\to\b_{1,0}$ in $W$. Finally,  we get a H-cycle $C$  for $\G$:
$$C:\quad \ph_{1,Y_4;0}\to\b_{1,0},\d_{Y_1,1}, \a\to\g_{Y_1,1},\ph_{1,Y_4;0}.$$\qed

\begin{figure}[htbp]
  \centering
  \begin{tikzpicture}[
    scale=0.75, transform shape,
    p/.style={draw, circle, fill=white, inner sep=0pt, minimum size=2.6mm},
    d/.style={fill=white, inner sep=2pt},
    box/.style={draw, rounded corners=6pt},
    lab/.style={anchor=north, yshift=3pt, font=\footnotesize},
    tlab/.style={anchor=south, font=\small},
    midarrow/.style={
      postaction=decorate,
      decoration={markings, mark=at position 0.5 with {\arrow{angle 45}}}
    },
    remidarrow/.style={
      postaction=decorate,
      decoration={markings, mark=at position 0.5 with {\arrowreversed{angle 45}}}
    }
  ]
    % 111
    \node[p] (a10)  at (0, 0)  {};
    \node[d] (d10)  at (2, 0)  {$\cdots$};
    \node[p] (a20)  at (4, 0)  {};
    \node[p] (a30)  at (6, 0)  {};
    \node[d] (d20)  at (8, 0)  {$\cdots$};
    \node[p] (a40)  at (10, 0)  {};
    \draw[remidarrow] (a10.east) -- (d10.west);
    \draw[remidarrow] (d10.east) -- (a20.west);
    %\draw[remidarrow] (a20.east) -- (a30.west);
    \draw[remidarrow] (a30.east) -- (d20.west);
    \draw[remidarrow] (d20.east) -- (a40.west);
    % 222
    \node[p] (a11)  at (0, -1)  {};
    \node[d] (d11)  at (2, -1)  {$\cdots$};
    \node[p] (a21)  at (4, -1)  {};
    \node[p] (a31)  at (6, -1)  {};
    \node[d] (d21)  at (8, -1)  {$\cdots$};
    \node[p] (a41)  at (10, -1)  {};
    \draw[remidarrow] (a11.east) -- (d11.west);
    \draw[remidarrow] (d11.east) -- (a21.west);
    \draw[remidarrow] (a21.east) -- (a31.west);
    \draw[remidarrow] (a31.east) -- (d21.west);
    \draw[remidarrow] (d21.east) -- (a41.west);
    % ...
    \node[d] (v1)  at (0, -2)  {$\vdots$};
    \node[d] (v2)  at (2, -2)  {$\vdots$};
    \node[d] (v3)  at (4, -2)  {$\vdots$};
    \node[d] (v4)  at (6, -2)  {$\vdots$};
    \node[d] (v5)  at (8, -2)  {$\vdots$};
    \node[d] (v6)  at (10, -2)  {$\vdots$};
    % 333
    \node[p] (a12)  at (0, -3)  {};
    \node[d] (d12)  at (2, -3)  {$\cdots$};
    \node[p] (a22)  at (4, -3)  {};
    \node[p] (a32)  at (6, -3)  {};
    \node[d] (d22)  at (8, -3)  {$\cdots$};
    \node[p] (a42)  at (10, -3)  {};
    \draw[remidarrow] (a12.east) -- (d12.west);
    \draw[remidarrow] (d12.east) -- (a22.west);
    \draw[remidarrow] (a22.east) -- (a32.west);
    \draw[remidarrow] (a32.east) -- (d22.west);
    \draw[remidarrow] (d22.east) -- (a42.west);

    % lian
    \draw[midarrow] (a10.east) -- (a41.west);
    \draw[midarrow] (a11.east) -- ($(a11)!0.36!(a42)$);
    \draw[midarrow] ($(a11)!0.64!(a42)$) -- (a42.west);
    \draw[midarrow, rounded corners=6pt]
    (a12.south) -- (0,-3.7) -- (10.9,-3.7) -- (10.9,0) -- (a40.east);

    % block
    \draw[box] (-0.7,0.5) rectangle (0.7,-3.5);
    \draw[box] (1.3,0.5) rectangle (2.7,-3.5);
    \draw[box] (3.3,0.5) rectangle (4.7,-3.5);
    \draw[box] (5.3,0.5) rectangle (6.7,-3.5);
    \draw[box] (7.3,0.5) rectangle (8.7,-3.5);
    \draw[box] (9.3,0.5) rectangle (10.7,-3.5);

    % sign
    \node[lab] at (a10.south) {$\d_{Y_1,1}$};
    \node[lab] at (a20.south) {$\d_{Y_0,1}$};
    \node[lab] at (a30.south) {$\g_{Y_0,1}$};
    \node[lab] at (a40.south) {$\g_{Y_2,1}$};
    \node[lab] at (a11.south) {$\d_{Y_1,2}$};
    \node[lab] at (a21.south) {$\d_{Y_0,2}$};
    \node[lab] at (a31.south) {$\g_{Y_0,2}$};
    \node[lab] at (a41.south) {$\g_{Y_2,2}$};
    \node[lab] at (a12.south) {$\d_{Y_1,0}$};
    \node[lab] at (a22.south) {$\d_{Y_0,0}$};
    \node[lab] at (a32.south) {$\g_{Y_0,0}$};
    \node[lab] at (a42.south) {$\g_{Y_2,0}$};

    \node[tlab] at (0,0.5) {$C_{Y_1}$};
    \node[tlab] at (4,0.5) {$C_{Y_0}$};
    \node[tlab] at (6,0.5) {$D_{Y_0}$};
    \node[tlab] at (10,0.5) {$D_{Y_1}$};

  \end{tikzpicture}
  \caption{{\footnotesize H-path $\d_{Y_0,1}\to\g_{Y_0,1}$ of $\G((\cup_Y C_Y)\cup (\cup_Y D_Y))$ in Lemma 6.4.(1).}}
%\end{figure}
\vskip 2mm
%\begin{figure}[htbp]
  \centering
  \begin{tikzpicture}[
    scale=0.75, transform shape,
    p/.style={draw, circle, fill=white, inner sep=0pt, minimum size=2.6mm},
    d/.style={fill=white, inner sep=2pt},
    box/.style={draw, rounded corners=6pt},
    lab/.style={anchor=north, yshift=3pt, font=\footnotesize},
    tlab/.style={anchor=south, font=\small},
    midarrow/.style={
      postaction=decorate,
      decoration={markings, mark=at position 0.5 with {\arrow{angle 45}}}
    },
    remidarrow/.style={
      postaction=decorate,
      decoration={markings, mark=at position 0.5 with {\arrowreversed{angle 45}}}
    }
  ]
    % 111
    \node[p] (a10)  at (0, 0)  {};
    \node[d] (d10)  at (2, 0)  {$\cdots$};
    \node[p] (a20)  at (4, 0)  {};
    \node[p] (a30)  at (6, 0)  {};
    \node[d] (d20)  at (8, 0)  {$\cdots$};
    \node[p] (a40)  at (10, 0)  {};
    \draw[remidarrow] (a10.east) -- (d10.west);
    \draw[remidarrow] (d10.east) -- (a20.west);
    %\draw[remidarrow] (a20.east) -- (a30.west);
    \draw[remidarrow] (a30.east) -- (d20.west);
    \draw[remidarrow] (d20.east) -- (a40.west);
    % 222
    \node[p] (a11)  at (0, -1)  {};
    \node[d] (d11)  at (2, -1)  {$\cdots$};
    \node[p] (a21)  at (4, -1)  {};
    \node[p] (a31)  at (6, -1)  {};
    \node[d] (d21)  at (8, -1)  {$\cdots$};
    \node[p] (a41)  at (10, -1)  {};
    \draw[remidarrow] (a11.east) -- (d11.west);
    \draw[remidarrow] (d11.east) -- (a21.west);
    \draw[remidarrow] (a21.east) -- (a31.west);
    \draw[remidarrow] (a31.east) -- (d21.west);
    \draw[remidarrow] (d21.east) -- (a41.west);
    % ...
    \node[d] (v1)  at (0, -2)  {$\vdots$};
    \node[d] (v2)  at (2, -2)  {$\vdots$};
    \node[d] (v3)  at (4, -2)  {$\vdots$};
    \node[d] (v4)  at (6, -2)  {$\vdots$};
    \node[d] (v5)  at (8, -2)  {$\vdots$};
    \node[d] (v6)  at (10, -2)  {$\vdots$};
    % 333
    \node[p] (a12)  at (0, -3)  {};
    \node[d] (d12)  at (2, -3)  {$\cdots$};
    \node[p] (a22)  at (4, -3)  {};
    \node[p] (a32)  at (6, -3)  {};
    \node[d] (d22)  at (8, -3)  {$\cdots$};
    \node[p] (a42)  at (10, -3)  {};
    \draw[remidarrow] (a12.east) -- (d12.west);
    \draw[remidarrow] (d12.east) -- (a22.west);
    \draw[remidarrow] (a22.east) -- (a32.west);
    \draw[remidarrow] (a32.east) -- (d22.west);
    \draw[remidarrow] (d22.east) -- (a42.west);

    % lian
    \draw[midarrow] (a10.east) -- (a41.west);
    \draw[midarrow] (a11.east) -- ($(a11)!0.36!(a42)$);
    \draw[midarrow] ($(a11)!0.64!(a42)$) -- (a42.west);
    \draw[midarrow, rounded corners=6pt]
    (a12.south) -- (0,-3.7) -- (10.9,-3.7) -- (10.9,0) -- (a40.east);

    % block
    \draw[box] (-0.7,0.5) rectangle (0.7,-3.5);
    \draw[box] (1.3,0.5) rectangle (2.7,-3.5);
    \draw[box] (3.3,0.5) rectangle (4.7,-3.5);
    \draw[box] (5.3,0.5) rectangle (6.7,-3.5);
    \draw[box] (7.3,0.5) rectangle (8.7,-3.5);
    \draw[box] (9.3,0.5) rectangle (10.7,-3.5);

    % sign
    \node[lab] at (a10.south) {$\ph_{1,Y_0;0}$};
    \node[lab] at (a20.south) {$\b_{2,0}$};
    \node[lab] at (a30.south) {$\b_{1,0}$};
    \node[lab] at (a40.south) {$\ph_{1,Y_1;0}$};
    \node[lab] at (a11.south) {$\ph_{1,Y_0;1}$};
    \node[lab] at (a21.south) {$\b_{2,1}$};
    \node[lab] at (a31.south) {$\b_{1,1}$};
    \node[lab] at (a41.south) {$\ph_{1,Y_1;1}$};
    \node[lab] at (a12.south) {$\ph_{1,Y_0;2^{2m}-2}$};
    \node[lab] at (a22.south) {$\b_{2,2^{2m}-2}$};
    \node[lab] at (a32.south) {$\b_{1,2^{2m}-2}$};
    \node[lab] at (a42.south) {$\ph_{1,Y_1;2^{2m}-2}$};

    \node[tlab] at (0,0.5) {$E_{1,Y_0}$};
    \node[tlab] at (4,0.5) {$B_2$};
    \node[tlab] at (6,0.5) {$B_1$};
    \node[tlab] at (10,0.5) {$E_{1,Y_1}$};

  \end{tikzpicture}
  \caption{{\footnotesize H-path $\b_{2,0}\to\b_{1,0}$ of $\G((\cup_Y E_{0,Y})\cup(\cup_Y E_{1,Y})\cup B_1\cup B_2)$ in Lemma 6.4.(1).}}
\end{figure}
\subsection{Action  on  maximal totally isotropic subspaces}

\begin{lem}\label{u2}
Let $M=\SU(4,q)$ or $\SU(5,q)$, acting on the set $\O$ of maximal totally isotropic subspaces. Then $\G(\O, \Del_s)$  has a H-cycle.
\end{lem}
\demo  (1) Let $M=\SU(4,q).$ Set $V=\lg u_1,u_2\rg\perp \lg v_1,v_2\rg$, where $\lg u_1,u_2\rg$ and $\lg v_1,v_2\rg$ are hyperbolic planes.
Then $|\O|=(q^3+1)(q+1)$.
For any $U\in \O$, we have
 $\D_s=\{W\in\O \di \dim(W\cap U)=1\}, \, {\rm  where}\, \,  |\D_s|=q(q^2+1).$
  Every isotropic point is contained in $q+1$ maximal totally isotropic subspaces. For an  isotropic point $\lg v\rg$, set $\o v=\{W\in\O \di v\in W\}$ so that  $|\o v|= q+1$.
 Let $u$ be an isotropic vector  in $\emph{}\lg u_1,u_2\rg$ and $\{\lg v_i\rg \di 0\le i\le q\}$  be all isotropic points in $\lg v_1,v_2\rg$. Assume $\d_{i,0}=v_i$ and  $\d_{i,j}=u+\th^j v_i$, where
 $0\le i \le q$ and $1\le j \le q^2-1$.
 The remaining arguments are completely same with that in   Lemma \ref{s2}.

\vskip 2mm
(2) Let $M=\SU(5,q).$ Set $V=\lg u_1,u_2\rg\perp \lg v_1,v_2\rg\perp \lg w\rg$, where $\lg u_1,u_2\rg$ and $\lg v_1,v_2\rg$ are hyperbolic planes and $f(w,w)=1$. Then $|\O|=(q^5+1)(q^3+1)$.
For any $U\in \O$,  we have
$\D_s=\{W\in\O \di \dim(W\cap U)=1\}, \, {\rm  where}\, \, |\D_s|=q^3(q^2+1).$
Note that every isotropic point is contained in $q^3+1$ maximal totally isotropic subspaces.
For an  isotropic point $\lg v\rg$, set $\o v=\{W\in\O \di v\in W\}$ so that  $|\o v|= q^3+1$.
Let $u$ be any isotropic vector in $\lg u_1,u_2\rg$, and $\{\lg v_i\rg \di 0\le i \le q^3\}$ be all isotropic points in $\lg v_1,v_2\rg\perp \lg w\rg$. Assume $\d_{i,0}=v_i$ and $\d_{i,j}=u+\th^j v_i$, where $0\le i \le q^3$ and $1\le j \le q^2-1$. Then we may  get a H-cycle for $\G$ by using the similar   arguments in  Lemma \ref{s2}.\qed

\end{document}